\RequirePackage[leqno]{amsmath}
\documentclass[leqno]{amsart}
\usepackage{amsmath,mathtools}
\makeatletter
\usepackage[utf8]{inputenc}     
\usepackage{multicol}
\usepackage[inline]{enumitem} 
\usepackage{color, graphics}
\usepackage[demo]{graphicx}
\usepackage{efbox,graphicx}
\efboxsetup{linecolor=gray,linewidth=0.8pt}
\usepackage{floatrow}
\usepackage{pgfplots}
\usepackage{latexsym, amsmath, amssymb}
\usepackage{caption}
\usepackage{subcaption}
\usepackage{pgfplots}
\usepackage{amsfonts}
\usepackage{pifont}
\usepackage{enumitem}
\setlist[itemize]{topsep=0pt,after=\vspace{1.5\baselineskip}}
\usepackage[colorlinks=true]{hyperref}
\usepackage{empheq}
\usepackage[T1]{fontenc}
\usepackage{tabularx}
\usepackage[leftcaption]{sidecap}
\sidecaptionvpos{figure}{m}
\usepackage{tikz}
\usepackage{caption}
\usepackage{tikz-cd}
\usetikzlibrary{decorations.markings}
\usetikzlibrary{decorations.pathmorphing}
\usepackage{tabularx,multirow,array} 
\usepackage[margin=0.5in]{geometry}
\usetikzlibrary{backgrounds,automata}

\usepackage{xparse}
\ExplSyntaxOn
\int_new:N \g_tohi_const_int
\int_new:N \g_tohi_const_sub_int
\tl_new:N  \g_tohi_const_char_tl

\cs_new_protected:Nn \tohi_print_constant:nn 
{
  #1 \textsubscript {#2} 
} 

\NewDocumentCommand\resetconstants{m}
{
 \int_gincr:N \g_tohi_const_int
 \int_gzero:N \g_tohi_const_sub_int
 \tl_gset:Nn  \g_tohi_const_char_tl {#1}
}

\NewDocumentCommand\const{m}
{
  \tl_if_exist:cTF
   {
    c_tohi_const_\int_use:N\g_tohi_const_int _#1_tl
   }
   {
    \tl_use:c {c_tohi_const_\int_use:N\g_tohi_const_int _#1_tl }
   }
   {
    \int_gincr:N \g_tohi_const_sub_int
    \tl_const:cx {c_tohi_const_\int_use:N\g_tohi_const_int _#1_tl }
     { \exp_not:N\tohi_print_constant:nn {\g_tohi_const_char_tl }{\int_use:N \g_tohi_const_sub_int}}
    \tl_use:c {c_tohi_const_\int_use:N\g_tohi_const_int _#1_tl }
   }
}

\ExplSyntaxOff

\newcommand{\inlineitem}[1][]{%
\ifnum\enit@type=\tw@
    {\descriptionlabel{#1}}
  \hspace{\labelsep}%
\else
  \ifnum\enit@type=\z@
       \refstepcounter{\@listctr}\fi
    \quad\@itemlabel\hspace{\labelsep}%
\fi}
\DeclarePairedDelimiter\abs{\lvert}{\rvert}
\DeclarePairedDelimiter\norm{\lVert}{\rVert}
\DeclarePairedDelimiter\tonda{(}{)}

\newcommand{\into}{\int_\Omega}

\setlist[itemize]{noitemsep, topsep=0pt}

\def\R{\mathbb R} \def\N{\mathbb N}

\def\R{\mathbb R} \def\N{\mathbb N} 
\def\TM{T_{max}} 
\def
\@cite
#1#2{[#1\if@tempswa, #2\fi]}
\makeatother
\newtheorem{theorem}{Theorem}[section]

\newtheorem{assumptions}[theorem]{Assumptions}

\newtheorem{lemma}[theorem]{Lemma}

\newtheorem{remark}{Remark}

\title[Blow-up solutions to Keller-Segel models] 
{Properties of given and detected unbounded solutions to a class of chemotaxis models}
\author[Alessandro Columbu, Silvia Frassu and Giuseppe Viglialoro]{}
\keywords{Chemotaxis, Attraction-repulsion, Nonlinear production, Blow-up time, Lower bound. \\
\textit{$^\star$Corresponding author}: silvia.frassu@unica.it}

\begin{document}

\maketitle

\centerline{\scshape{\scshape{Alessandro Columbu, Silvia Frassu$^{\star}$ \and Giuseppe Viglialoro}}}
\medskip

{
\medskip
\centerline{Dipartimento di Matematica e Informatica}
\centerline{Universit\`{a} di Cagliari}
\centerline{Via Ospedale 72, 09124. Cagliari (Italy)}
\medskip
}
\bigskip
\begin{abstract}
This paper deals with unbounded solutions to a class of chemotaxis systems. In particular, for a rather general attraction-repulsion model, with nonlinear productions, diffusion, sensitivities and logistic term, we detect Lebesgue spaces where given unbounded solutions blow-up also in the corresponding norms of those spaces; subsequently, estimates for the blow-up time are established. Finally, for a simplified version of the model, some blow-up criteria are proved.

More precisely, we analyze a zero-flux chemotaxis system essentially described as 
\begin{equation}\label{problem_abstract}
\tag{$\Diamond$}
\begin{dcases}
u_t= \nabla \cdot ( (u+1)^{m_1-1}\nabla u -\chi u(u+1)^{m_2-1}\nabla v + \xi u(u+1)^{m_3-1}\nabla w) +\lambda u -\mu u^k  & \text{ in } \Omega \times (0,\TM),\\
0= \Delta v -\frac{1}{\abs{\Omega}}\into  u^\alpha   + u^\alpha = \Delta w - \frac{1}{\abs{\Omega}}\into  u^\beta  + u^\beta & \text{ in } \Omega \times (0,\TM).
\end{dcases}
\end{equation}
The problem is formulated in a bounded and smooth domain $\Omega$ of $\R^n$, with $n\geq 1$, for some $m_1,m_2,m_3\in \R$, $\chi, \xi, \alpha,\beta, \lambda,\mu>0$, $k >1$, and with $\TM\in (0,\infty]$. A sufficiently regular initial data $u_0\geq 0$ is also fixed.  

Under specific relations involving the above parameters, one of these  always requiring some largeness conditions on $m_2+\alpha$,    
\begin{enumerate}[label=\roman*)]
\item \label{Abs:Obj-1} we prove that any given solution to \eqref{problem_abstract}, blowing-up at some finite time $\TM$, becomes unbounded also in $L^{\mathfrak{p}}(\Omega)$-norm, for all ${\mathfrak{p}}>\frac{n}{2}\left(m_2-m_1+\alpha\right)$; 
\item \label{Abs:Obj-2} we give lower bounds $T$ (depending on $\int_\Omega u_0^{\Bar{p}}$) of $\TM$ for the aforementioned solutions in some $L^{\Bar{p}}(\Omega)$-norm, being $\Bar{p}=\Bar{p}(n,m_1,m_2,m_3,\alpha,\beta)\geq\mathfrak{p}$;    
\item \label{Abs:Obj-3} whenever  $m_2=m_3$, we establish sufficient conditions on the parameters ensuring that for some $u_0$ solutions to \eqref{problem_abstract} effectively are unbounded at some finite time. 
\end{enumerate}  
Within the context of blow-up phenomena connected to problem \eqref{problem_abstract}, this research partially improves  the analysis in \cite{wang2023blow} and moreover contributes to enrich the level of knowledge on the topic. 
\end{abstract}
\resetconstants{c}
\section{Introduction, motivations and state of the art}\label{Intro}
\subsection{The continuity equation; the initial-boundary value problem}
The well known continuity equation 
\begin{equation}\label{ContinuityEq}
    u_t=\nabla \cdot F+h
\end{equation}
describes the transport of some quantity $u=u(x,t)$, at the position $x$ and at the time $t>0$. In this equation the flux $F$ models the motion of such a quantity, whereas $h$ is an additional source idealizing some external action by means of which $u$ itself may be created or destroyed throughout the time.

In this paper we are interested in the analysis of equation \eqref{ContinuityEq} in the context of self organization mechanisms for biological populations, i.e. phenomena for which organisms or entities direct their trajectory in response to one or more chemical stimuli.  More precisely, we want to deal with the motion of a certain cell density $u=u(x, t)$ whose flux has a smooth diffusive part and another contrasting this spread. In the specific, this counterpart accounts of an attractive and a repulsive effect, associated to two chemical signals and indicated respectively with $v=v(x,t)$ and $w=w(x,t)$; $v$ (the \textit{chemoattractant}) tends to gather the cells, $w$ (the \textit{chemorepellent}) to scatter them. Additionally an external source with an increasing and decreasing effect on the cell density is also included.  

For our purposes, wanting to formulate what said above in terms of  the continuity equation \eqref{ContinuityEq}, it appears meaningful (and convenient) defining for $\chi,\xi,\lambda,\mu>0$, $m_1,m_2,m_3\in \R$ and $k>1$ the fluxes $F=F_{m_1,m_2,m_3}=F_{m_1,m_2,m_3}(u,v,w)$, $G_{m_1}=G_{m_1}(u)$, $H_{m_2}=H_{m_2}(u,v)$, $I_{m_3}=I_{m_3}(u,w)$ and the source $h=h_k=h_k(u)$:
\begin{equation}\label{FluxAndSource} 
\begin{cases}
F=F_{m_1,m_2,m_3}=(u+1)^{m_1-1}\nabla u -\chi u(u+1)^{m_2-1}\nabla v+\xi u(u + 1)^{m_3-1}\nabla w=:G_{m_1} +H_{m_2}+I_{m_3}\\ 
  h= h_k=\lambda u -\mu u^k.
    \end{cases}
\end{equation}
 In this way, the diffusion is smoother and smoother for higher and higher values of $m_1$, the aggregation/repulsion effects 
$-\chi u(u+1)^{m_2-1}\nabla v/\xi u(u + 1)^{m_3-1}\nabla w$ increase for larger sizes of $\chi$ and $\xi$ and $m_2$ and $m_3$, and the cell density may increment with rate $\lambda u$ and may  attenuate with rate  $-\mu u^k$. Naturally, since the flux is influenced by the two signals $v$ and $w$, we will have to consider two more equations (which for the time being we indicate with $P(v)=0$ and $Q(w)=0$, but that will be specified later), one for the chemoattractant $v$ and another for the chemorepellent $w$, to be coupled with the continuity equation. Furthermore, the analysis is studied in impenetrable domains (so homogeneous Neumann  or zero flux boundary conditions are imposed) and some initial configurations for the cell and chemical densities are assigned: essentially, with position \eqref{FluxAndSource} in mind, we are concerned with this initial boundary value problem: 
\begin{equation}\label{problemAstratto}
\begin{cases}
        u_t=\nabla \cdot F_{m_1,m_2,m_3}+h_k & \textrm{ in } \Omega \times (0,T_{max}), \\
        P(v)=Q(w)=0 & \textrm{ in } \Omega \times (0,T_{max}),\\
        u_0(x)=u(x,0)\geq 0; \;v_0(x)=v(x,0)\geq 0;\; w_0(x)=w(x,0) \geq0 & x \in \bar\Omega,\,\\
        u_{\nu}=v_{\nu}=w_{\nu}=0 & \textrm{ on } \partial \Omega \times (0,T_{max}).
\end{cases}
\end{equation}
The problem is formulated in a bounded and smooth domain $\Omega$ of $\R^n$, with $n\geq 1$, $u_\nu$ (and similarly for $v_\nu$ and $w_\nu$) indicates the outward normal derivative of $u$ on $\partial \Omega$. Moreover,  $\TM\in (0,\infty]$  identifies the maximum time up to which solutions to the system can be extended. 
 \subsection{A view on the state of the art: the attractive and the repulsive models and the attraction-repulsion model} 
The aforementioned discussion finds, of course, its roots in the well-known Keller--Segel models idealizing chemotaxis phenomena (see the celebrated papers \cite{K-S-1970,Keller-1971-MC,Keller-1971-TBC}), that since the last 50 years have been attracting the interest of the mathematical community.

In particular, if we refer to chemotaxis models with single proliferation signal, taking in mind  \eqref{FluxAndSource},  problem \eqref{problemAstratto} is a (more general) combination of this aggregative signal-production mechanism 
\begin{equation}\label{problemOriginalKS} 
u_t= \nabla \cdot \left(G_1+H_1\right)= \Delta u- \chi \nabla \cdot (u \nabla v)\quad \textrm{and} \quad 
P(v)=P^\tau_1(v)=\tau v_t-\Delta v+v-u=0,  \quad  \textrm{ in } \Omega \times (0,T_{max}),
\end{equation}
and this repulsive signal-production one 
\begin{equation}\label{problemOriginalKSCosnumption}
u_t=  \nabla \cdot \left(G_1+I_1\right)= \Delta u+ \xi \nabla \cdot (u \nabla w)  \quad \textrm{and} \quad 
Q(w)=Q^\tau_1(w)=\tau w_t-\Delta w+w -u=0,  \quad  \textrm{ in } \Omega \times (0,T_{max}).
\end{equation}
(Here $\tau\in\{0,1\}$ and it distinguishes between a stationary and evolutive equation for the chemical.)
The above models present linear diffusion and linear production rates; specifically,  $v$ and $w$ are linearly produced by the cells themselves, and their mechanism is opposite when in $P^\tau_1(v)$ the term $v-u$ (or in $Q_1(w)$ the term $w-u$) is replaced by $uv$ (or $uw$); in this case the particle density consumes the chemical. (We will spend only few words on models with absorption.) As far as problem \eqref{problemOriginalKS} is concerned, since the attractive signal $v$ increases with $u$, the natural spreading process of the cells' density could interrupt and very high and spatially concentrated spikes formations (\textit{chemotactic collapse} or \textit{blow-up at finite time}) may appear; this is, generally, due to the size of the chemosensitvity $\chi$, the initial mass of the particle distribution, i.e.,  $m=\int_\Omega u_0(x)dx,$ and the space dimension $n$. In this direction, the reader interested in learning more can find in \cite{HerreroVelazquez,JaLu,Nagai,WinklAggre} analyses dealing with existence and properties of global, uniformly bounded, or blow-up (local) solutions to models connected to \eqref{problemOriginalKS}.

On the other hand, for nonlinear segregation chemotaxis models like those we are interested in, when in problem \eqref{problemOriginalKS} one has that $P(v)=P^1_\alpha(v)=v_t-\Delta v+v-u^\alpha=0$, with  $0<\alpha<\frac{2}{n}$ ($n\geq1$), uniform boundedness of all its solutions is proved in \cite{LiuTaoFullyParNonlinearProd}.

Concerning the literature about   problem \eqref{problemOriginalKSCosnumption}, it seems  rather poor and general (see, for instance, \cite{Mock74SIAM,Mock75JMAA} for analyses on similar contexts). In particular, no result on the blow-up scenario is available; this is meaningful due to the repulsive nature of the phenomenon.

Contrarily, the level of understanding for attraction-repulsion chemotaxis problems involving both \eqref{problemOriginalKS} and \eqref{problemOriginalKSCosnumption} is sensitively rich; more specifically, if we refer to the linear diffusion and sensitivities version of model \eqref{problemAstratto}, for which $F=F_{1,1,1}$ and $P^\tau_\alpha(v)=\tau v_t-\Delta v+b v -a u^\alpha$ and $Q^\tau_\beta(w)=\tau w_t- \Delta w + d w - c u^\beta$,
$a,b,c,d,\alpha,\beta>0$, equipped with regular initial data $u_0(x),\tau v_0(x),\tau w_0(x)\geq 0$, we can recollect the following outcomes. In the absence of logistics ($h_k\equiv 0$), when linear growths of the chemoattractant and the chemorepellent are taken into consideration, and for elliptic equations for the chemicals (i.e. when $P^0_1(v)=Q^0_1(w)=0$), the value $\Theta:=\chi a-\xi c$ measures the difference between the attraction and repulsion impacts, and it is such that whenever $\Theta<0$ (repulsion-dominated regime), in any dimension all solutions to the model are globally bounded, whereas for $\Theta>0$ (attraction-dominated regime) and $n=2$ unbounded solutions can be detected (see \cite{GuoJiangZhengAttr-Rep,LI-LiAttrRepuls,TaoWanM3ASAttrRep,VIGLIALORO-JMAA-BlowUp-Attr-Rep,YUGUOZHENG-Attr-Repul} for some details on the issue). Indeed, for more general expressions of the proliferation laws, modelled by the  equations $P^0_\alpha(v)=Q^0_\beta(w)=0$, to the best of our knowledge, \cite{ViglialoroMatNacAttr-Repul} is  the most recent result in this direction; herein some interplay between $\alpha$ and $\beta$ and some technical conditions on $\xi$ and $u_0$ are established so to ensure globality and boundedness of classical solutions.  (See also \cite{ChiyoYokotaBlow-UpAttRe} for blow-up results in the frame of nonlinear attraction-repulsion models with logistics as those formulated in \eqref{problemAstratto} with $F=F_{m_1,m_2,m_3}$ and $h=h_k$, and with linear segregation for the stimuli, i.e.  with equations for $v,w$  reading as $P_1^0(v)=Q_1^0(w)=0$.)

Putting our attention on evolutive equations for chemoattractant and chemorepellent, $P^1_1(v)=Q^1_1(w)=0$, in \cite{TaoWanM3ASAttrRep} it is proved that in two-dimensional domains sufficiently smooth initial data emanate global-in-time bounded solutions whenever  
$$\Theta<0 \quad \textrm{and} \quad b=d  \quad \textrm{or} \quad  \Theta<0 \quad \textrm{and} \quad  -\frac{\chi^2\alpha^2(b-d)^2}{2\Theta {b}^2 C} \int_\Omega u_0(x)dx\leq 1, \quad \textrm{for some }\; C>0.$$
(As to blow-up results we are only aware of \cite{LankeitJMAA-BlowupParabolicoAttr-Rep}, where unbounded solutions in three-dimensional domains are constructed.) 
When $h=h_k\not\equiv 0$, for both linear and nonlinear productions scenarios, and stationary or evolutive equations (formally, $P^{\tau}_\alpha(v)=Q^{\tau}_\beta(w)=0$),  criteria toward boundedness, long time behaviors and  blow-up issues for related solutions are studied in \cite{LiangEtAlAtt-RepNonLinProdLogist-2020,XinluEtAl2022-Asymp-AttRepNonlinProd,ChiyoMarrasTanakaYokota2021,GuoqiangBin-2022-3DAttRep}.
\subsection{The nonlocal case} In this work we are mainly interested to the so-called nonlocal models tied to \eqref{problemAstratto}, for which, specifically,  $P(v)=P_\alpha(v):=\Delta v-\frac{1}{|\Omega|}\int_\Omega u^\alpha +u^\alpha=0$,  and $Q(w)=Q_\beta(w):=\Delta w-\frac{1}{|\Omega|}\int_\Omega u^\beta +u^\beta$. In particular, recalling the position in \eqref{FluxAndSource},  we herein mention the most recent researches we are aware about and inspiring our study; to this purpose we refer to the only attraction  version
\begin{equation}\label{OnlyATTractionproblem} 
u_t= \nabla \cdot \left(G_{m_1}+H_{m_2}\right) +h_k \quad \textrm{and} \quad 
P_\alpha(v)=0  \quad  \textrm{ in } \Omega \times (0,T_{max}),
\end{equation}
and the attraction-repulsion one  
\begin{equation}\label{ATTract-Repulsionproblem}  
u_t= \nabla \cdot F_{m_1,m_2,m_3} +h_k  \quad \textrm{and} \quad 
P_\alpha(v)=Q_\beta(w)=0  \quad  \textrm{ in } \Omega \times (0,T_{max}).
\end{equation}
For the linear diffusion and sensitivity version of model \eqref{OnlyATTractionproblem}, i.e. with flux $G_{1}+H_{1}$, if $h_k\equiv 0$, it is known that  boundedness of solutions is achieved for any $n\geq 1$ and $0<\alpha<\frac{2}{n}$, whereas for $\alpha>\frac{2}{n}$ blow-up phenomena may be observed (see \cite{WinklerNoNLinearanalysisSublinearProduction}). When dampening logistic terms take part in the mechanism, for the linear flux $G_1+H_1$ appearances of $\delta$-formations at finite time have been detected for some sub-quadratic growth of $h_k$, and precisely for $h_k$ with $1<k<\frac{n(\alpha+1)}{n+2}<2$ (see \cite{yi2021blow}). But there is more; for the limit linear production scenario, corresponding to $P_1(v)=0$, some unbounded solutions have been constructed in \cite{FuestCriticalNoDEA} even for quadratic sources $h=h_2$, whenever $n\geq 5$ and $\mu\in \left(0,\frac{n-4}{n}\right).$ 

In nonlinear models without dampening logistic effects (i.e., general flux  $G_{m_1}+H_{m_2}$ and $h_k\equiv 0$) and linear production (i.e. $P_1(v)=0$), in \cite{winkler2010boundedness} it is shown, among other things, that  for $m_1\leq 1$, $m_2>0$, $m_2>m_1+\frac{2}{n}-1$ situations with unbounded solutions at some finite time $T_{max}$ can be found. (See also \cite{MarrasNishinoViglialoro} for questions connected to estimates of $\TM$.) Some results have been also extended in \cite{TANAKA2021} when $h_k\not\equiv 0$ and for nonlinear segregation contexts, $P_\alpha(v)=0$; in particular, inter alia, for $m_1\in \R$, $m_2>0$,  blow-up phenomena are seen to appear if $ m_2+\alpha > \max\bigl\{m_1+\frac{2}{n}k, k\bigr\}$, whenever  $m_1 \geq 0$ or $
m_2+\alpha > \max\bigl\{\frac{2}{n}k, k\bigr\}$, provided  $m_1 <0$. 

On the other hand, for the attraction-repulsion models, in \cite{LiuLiBlowUpAttr-Rep} it is proved, together with other results, that if $P_\alpha(v)=Q_\beta(w)=0$, $\alpha>\frac{2}{n}$ and $\alpha>\beta$ ensure the existence of unbounded solutions to \eqref{ATTract-Repulsionproblem} for the  linear flux   $F=F_{1,1,1}$, without logistic ($h_k\equiv 0$). Conversely, detecting gathering mechanisms for the nonlinear situation is more complex and we are only aware of \cite{wang2023blow}; indeed this issue is therein addressed only for $F=F_{m_1,1,1}$, with $m_1\in \R$, but even in presence of dampening logistics. Since in our research we will show the existence of blow-up solutions for a larger class of fluxes, precisely for $F=F_{m_1,m_2,m_2}$ with $m_1\in \R$ and any $m_2=m_3>0$, we will spend more words to analyze details of \cite{wang2023blow} below, precisely in $\S$\ref{SubsectionImproving}.

\section{Presentation of the model and of the main results. Aims of the paper}
\subsection{The formulation of the mathematical problem}
In light of what we presented so far, and under the aforementioned main positions,  basically in this paper we are interested in properties of unbounded classical solutions $(u,v,w)=(u(x,t),v(x,t),w(x,t))$ to this problem 
\begin{equation}\label{problem3x3}
\begin{dcases}
u_t= \nabla \cdot ( (u+1)^{m_1-1}\nabla u -\chi u(u+1)^{m_2-1}\nabla v + \xi u(u+1)^{m_3-1}\nabla w) +\lambda u -\mu u^k & \text{ in } \Omega \times (0,\TM),\\
0= \Delta v - m_1(t)  + f_1(u) & \text{ in } \Omega \times (0,\TM),\\
0= \Delta w - m_2(t)  + f_2(u) & \text{ in } \Omega \times (0,\TM),\\
u_{\nu}=v_{\nu}=w_{\nu}=0 & \text{ on } \partial \Omega \times (0,\TM),\\
u(x,0)=u_0(x) & x \in \bar\Omega,\\
\int_\Omega v(x,t)dx=\int_\Omega w(x,t)dx=0 & \text{ for all } t\in  (0,\TM).
\end{dcases}
\end{equation}
 Additionally, the initial cell distribution $u_0:=u_0(x)$ and the production laws $f_i=f_i(u)$ (for $i \in \{1,2\}$) are supposed to be nonnegative and sufficiently regular and  the functions $m_i(t)$ are defined for compatibility in the second and third equation (by integrating these over $\Omega$) in terms of $f_i(u)$ themselves, exactly as
\begin{equation}\label{Definitionm1andme}
   m_1(t)=\frac{1}{\abs{\Omega}}\into f_1(u)\quad \textrm{and} \quad m_2(t)=\frac{1}{\abs{\Omega}}\into f_2(u).
\end{equation}
Herein we assume 
\begin{equation}\label{f_i_u0_GENERALE}
0\leq f_i \in \bigcup_{\theta \in (0,1)} C^{\theta}_{loc}\tonda*{[0,\infty)} \cap C^1((0,\infty)) \;\; \text{ and } 0\leq u_0 \in \bigcup_{\theta \in (0,1)} C^{\theta}(\bar{\Omega});
\end{equation} 
 we also might need that for all $s\geq 0$, $\alpha, \beta>0$ and some $k_1, k_2,k_3>0$,
\begin{equation}\label{CondizioniSuFAlessandro}
f_1(s)\leq k_1(s+1)^\alpha \quad \textrm{and} \quad f_2(s)\leq k_2(s+1)^\beta, \quad 
\end{equation}
or 
\begin{equation}\label{f_iDef}
f_1,f_2 \,\text{ nondecreasing}, f_1(s) \geq k_3 (s+1)^{\alpha},  \quad  f_2(s) \leq k_2 (s+1)^{\beta} \textrm{ and $u_0=u_0(|x|)$ radially symmetric and nonincreasing}. 
\end{equation}
\begin{remark}
Let us clarify that in nonlocal models, and henceforth in this work, $v$ stands for the deviation of the chemoattractant; the deviation is  the difference between the signal concentration and its mean, and that it changes sign in contrast to what happens with the cell and signal densities (which are nonnegative). In particular, it follows from the definition of $v$ itself that its mean is zero (as specified in the last positions of the problem \eqref{problem3x3}), which in turn ensures the uniqueness of the solution of the Poisson equation under homogeneous Neumann boundary conditions. The same comments apply for the chemorepellent $w$. (We did not introduce different symbols to indicate the chemicals and their deviations since  it is clear from the context.)
\end{remark}
\subsection{Presentation of the Theorems. Overall aims of the paper} 
Our project finds its motivations in the observation that there is no automatic connection between the occurrence of blow-up for solutions to model \eqref{problem3x3} in the  $L^\infty(\Omega)$-norm and that in $L^p(\Omega)$-norm ($p>1$). Indeed, for a bounded domain $\Omega$, it is seen that  \[\| u(\cdot,t)\|_{L^p(\Omega)}\leq |\Omega|^\frac{1}{p}\|u(\cdot,t)\|_{L^\infty(\Omega)},\]
so that unboundedness in $L^p(\Omega)$-norm implies that in $L^\infty(\Omega)$-norm, but oppositely $\int_\Omega u^p$ might even remain bounded in a neighborhood of $T_{max}$ when $\max_{\Omega}u$ uncontrollably increases at some finite time $T_{max}$.  

In light of this, in order to bridge the gap between the analysis of the blow-up time $T_{max}$ in the two different mentioned norms, we aim at 
\begin{enumerate}[label=(\roman*)]
\item \label{Aim1}   detecting suitable $L^p(\Omega)$ spaces, for certain $p$ depending on $n,m_1,m_2,m_3,\alpha$ and $\beta$, such that given unbounded solutions also blow up in the associated  $L^p(\Omega)$-norms; 
\item \label{Aim2}  providing lower bounds for the blow-up time of the aforementioned solutions in these $L^p(\Omega)$-norms.
\end{enumerate}
Additionally, another objective of our work is
\begin{enumerate}[label=(\roman*)]
 \setcounter{enumi}{2}
    \item \label{Aim3} giving sufficient conditions on the data of the model  such that related solutions are actually unbounded at finite time.
\end{enumerate}
In order to deal with issues \ref{Aim1}, \ref{Aim2} and \ref{Aim3}, we fix the following relations, determining some precise  interplay involving  constants defining system \eqref{problem3x3}: 
\begin{assumptions} \label{generalassumptions}
Let $n\in\N$, $m_1,m_2,m_3\in\R$ and $\alpha,\beta,\chi,\xi,\lambda,\mu>0$ and $k>1$ be such that 
\begin{enumerate}[label={($\mathcal{H}_{\arabic*}$)}]
    \item \label{condizioneprincipale} $m_2+\alpha>m_1+\frac 2n$.
\end{enumerate}
Moreover, for $\beta>0$, let either
\begin{enumerate}[resume, label={($\mathcal{H}_{\arabic*}$)}]
        \item \label{As1}  $m_2+\alpha>\max\{1,m_3+\beta\}$ \quad or 
        \inlineitem \label{As2} $m_2+\alpha\geq m_3+\beta$ and $m_1>1-\frac 2n$,
        \end{enumerate}
whereas, for $\beta\in(0,1]$, let either
\begin{enumerate}[resume, label={($\mathcal{H}_{\arabic*}$)}]
        \item \label{As3} $m_3\leq 1$ and $m_1>1-\frac 2n$ \quad \,\, or
        \inlineitem \label{As4} $m_2+\alpha \geq  m_3$ and $m_1>1-\frac 2n$.
    \end{enumerate}
Finally, for $m_2, \beta>0$, let also
\begin{enumerate}[resume, label={($\mathcal{H}_{\arabic*}$)}]
\item \label{BU}
$\alpha>\beta \quad \text{and} \quad 
\begin{cases}
m_2+\alpha > \max\bigl\{m_1+\frac{2}{n}k, k\bigr\} & \text{ if } m_1 \geq 0,\\
m_2+\alpha > \max\bigl\{\frac{2}{n}k, k\bigr\} & \text{ if } m_1 < 0.
\end{cases}
$
\end{enumerate}
\end{assumptions}
With the support of the above position, and as far as the analysis of \ref{Aim1} is concerned, in the spirit of \cite[Theorem 2.2]{FREITAGLpLinftyCoincide}
we will prove this first result, dealing with properties of \textit{given unbounded solutions} to model \eqref{problem3x3}. Specifically, the proof is based on the analysis of the functional $\varphi(t)=\frac{1}{p}\int_\Omega (u+1)^p$, defined for local solutions to problem \eqref{problem3x3} on $(0,\TM)$. We will show that if for some $p$ sufficiently large $\varphi(t)$ is uniformly bounded in time, then it also is for any arbitrarily large  $p>1$, so contrasting with the unboundedness of $u$ itself.
\begin{theorem}\label{maintheorem1} Let $\Omega$ be a bounded and smooth domain of $\R^n$,  and  condition \ref{condizioneprincipale} as well as one of \ref{As1},\ref{As2},\ref{As3}, \ref{As4} in Assumptions \ref{generalassumptions} hold true. Moreover, for $f_i$ and $u_0$ complying with \eqref{f_i_u0_GENERALE} and \eqref{CondizioniSuFAlessandro},  let $(u,v,w)$, with 
\[
u \in C^0(\bar{\Omega} \times [0,\TM)) \cap C^{2,1}(\bar{\Omega} \times (0,\TM)) \;\textrm{and}\;
v, w \in \displaystyle \bigcap_{q>n} L^{\infty}((0,\TM); W^{1,q}(\Omega)) \cap C^{2,0}(\bar{\Omega} \times (0,\TM)), 
\]
be a solution to problem \eqref{problem3x3} 
which blows-up at some finite time $\TM$, in the sense that 
\begin{equation}\label{blowup}
    \limsup_{t\to \TM} {\lVert u(\cdot,t)\rVert_{L^\infty(\Omega)}}=+\infty.
\end{equation}
Then, for any $\mathfrak{p}>\frac n2 \left(m_2-m_1+\alpha\right)$, we also have 
\begin{equation*}
    \limsup_{t\to \TM} {\lVert u(\cdot,t)\rVert_{L^\mathfrak{p}(\Omega)}}=+\infty.
\end{equation*}
\end{theorem}
Successively aim \ref{Aim2} is achieved by establishing for the same functional $\varphi(t)$ a first order differential inequality
(ODI) of the type $\varphi'(t)\leq \Psi(\varphi(t))$ on $(0,T_{max})$. In particular, for any $\tau_0>0$ the function $\Psi(\tau)$
obeys the Osgood criterion (\cite{Osgood}),
\begin{equation}\label{OsgoodCriterion}
\int_{\tau_0}^\infty \frac{d\tau}{\Psi(\tau)}<\infty\quad \textrm{with } \tau_0>0,
\end{equation}
so that an integration on $(0,T_{max})$ of the ODI 
implies, whenever $\limsup_{t\rightarrow \TM} \varphi(t)=\infty$, the following lower bound for $\TM$ $$T_{max}\geq \int_{\varphi(0)}^\infty \frac{d \varphi}{\Psi(\varphi)}:=T,$$ and thereby a safe interval of existence $[0, T)$ for solutions to the model itself. 
\begin{theorem}\label{maintheorem2}
Let the hypotheses of Theorem \ref{maintheorem1} be satisfied. Then there exist $\Bar{p}>1$ and positive constants $\mathcal{A},\mathcal{B},\mathcal{C}, \mathcal{D}$, as well as $\gamma>\delta>1$, such that the blow-up time $\TM$ complies with both the implicit estimate
\begin{equation}\label{stimatempobu}
    \TM\geq\int_{\varphi(0)}^{\infty} \frac{d\tau}{\mathcal{A}\tau^\gamma+\mathcal{B}\tau^\delta+\mathcal{C}},
\end{equation}
and the explicit one
\begin{equation*}
\label{stimatempobuExplicit}
    \TM\geq \frac{\varphi(0)^{1-\gamma}}{\mathcal{D}(\gamma-1)},
\end{equation*}
where $\varphi(0)=\displaystyle\frac 1p \into (u_0+1)^{p}$, for any $p\geq\Bar{p}$.
\end{theorem}
Finally, the last theorem (connected to item \ref{Aim3}) establishes (at least in a particular case) the \textit{existence of unbounded solutions} to system \eqref{problem3x3},  so as to make the two previous statements meaningful. The basic idea consists on the analysis of the temporal evolution of the functional $\phi(t):=\int_0^{s_0} s^{-\gamma} (s_0-s) U(s,t)\,ds$ for  $t \in [0,\TM)$, being $U(s,t)$ the so-called mass accumulation function of $u$,  obeying a superlinear ODI.  
\begin{theorem}\label{ThBlowUp}
Let $\chi,\xi,\lambda,\mu,m_2=m_3>0$ and $k>1$. Additionally, for some $R>0$, let $\Omega=B_R(0) \subset \R^n$ be a ball, $u_0$, $f_1$ and $f_2$ satisfy \eqref{f_i_u0_GENERALE} and \eqref{f_iDef}. Finally, let hypotheses \ref{BU} be satisfied. Then for any $M_0>C$, with $C=\left(\frac{\lambda}{\mu} |\Omega|^{k-1}\right)^{\frac{1}{k-1}}$, there exist $\epsilon_0 \in (0,M_0)$ and $r_*\in(0,R)$ with the property that whenever $u_0$ is also chosen such  that 
\begin{equation}\label{M0}
\int_\Omega u_0(x)dx=M_0 \quad \text{and} \quad \int_{B_{r_*}(0)} u_0(x) dx\geq M_0-\epsilon_0,
\end{equation}
the corresponding classical solution $(u,v,w)$ to model \eqref{problem3x3} blows up at some finite time $\TM$, in the sense of relation \eqref{blowup}.
\end{theorem} 
\section{Miscellaneous and general comments}\label{section3}
\subsection{On the parameters $\mathfrak{p}$ and $\bar{p}$} We herein want to spend some words on the role of the parameters $\mathfrak{p}$ and $\bar{p}$ appearing in Theorems \ref{maintheorem1} and \ref{maintheorem2}. In particular, as we will observe in Lemma \ref{Lemmapbarra}   below, $\mathfrak{p}$ and $\bar{p}$ depend on $n, m_1,m_2,m_3,\alpha$ and $\beta$, and Figure \ref{prova5} collects some of their values on the $p$ axis. In the specific, for a blowing-up solution to \eqref{problem3x3}, we will have that $\int_\Omega u^p$ is bounded on $(0,\TM)$ for $p=1$, while it blows up for $p\geq \mathfrak{p}$; in general, the behaviour of $\int_\Omega u^p$ on $(0,\TM)$ for $p\in\left(1,\mathfrak{p}\right)$ is unknown. On the other hand, an estimate for $\TM$ is given in terms of $\frac{1}{p}\int_\Omega \left(u_0+1\right)^p$, for $p\geq \Bar{p}$, but not when $p\in\left(\mathfrak{p},\Bar{p}\right)$. 
\begin{figure}[h!]
    \centering 
    \begin{tikzpicture}
        \filldraw[black](0,0)circle (1.8pt)--(0.37,0);
        \draw[dashed](0.37,0) -- (2.7,0)node{$\circ$}  --(2.7,0);
        \draw[decorate, decoration={snake, segment length=1mm, amplitude=0.25mm}](2.7,0) -- (14.99,0);
        \draw[->](14.99,0) -- (15,0)node[right]{$ p $};
        \draw(0.37,0)--(0.37,0)node{$\circ$};
        \node[pin={[pin distance=0.2cm, pin edge={black}]{$1$}}] at (0,0) {};
        \node[pin={[pin distance=0.2cm, pin edge={black}]{$\mathfrak{p}$}}] at (0.37,0) {};
        \node[pin={[pin distance=0.2cm, pin edge={black}]{$\Bar{p}$}}] at (2.7,0) {};
    \end{tikzpicture}
    \vskip 0.25cm
    \begin{tikzpicture}
        \filldraw[black](0,0)circle (1.8pt)--(2,0);
        \draw[dashed](2,0) -- (13.7,0)node{$\circ$}  --(13.7,0);
        \draw[decorate, decoration={snake, segment length=1mm, amplitude=0.25mm}](13.7,0) -- (14.99,0);
        \draw[->](14.99,0) -- (15,0)node[right]{$ p $};
        \draw(2,0)node{$\circ$};
        \node[pin={[pin distance=0.2cm, pin edge={black}]{$1$}}] at (0,0) {};
        \node[pin={[pin distance=0.2cm, pin edge={black}]{$\mathfrak{p}$}}] at (2,0) {};
        \node[pin={[pin distance=0.2cm, pin edge={black}]{$\Bar{p}$}}] at (13.7,0) {};
    \end{tikzpicture}
    \vskip 0.25cm
    \begin{tikzpicture}
        \filldraw[black](0,0)circle (1.8pt)--(2.47,0);
        \draw[dashed](2.47,0) -- (3.85,0)node{$\circ$}  --(3.85,0);
        \draw[decorate, decoration={snake, segment length=1mm, amplitude=0.25mm}](3.85,0) -- (14.99,0);
        \draw[->](14.99,0) -- (15,0)node[right]{$ p $};
        \draw(2.47,0)node{$\circ$}; 
        \node[pin={[pin distance=0.2cm, pin edge={black}]{$1$}}] at (0,0) {};
        \node[pin={[pin distance=0.2cm, pin edge={black}]{$\mathfrak{p}$}}] at (2.47,0) {};
        \node[pin={[pin distance=0.2cm, pin edge={black}]{$\Bar{p}$}}] at (3.85,0) {};
    \end{tikzpicture}
    \vskip 0.25cm
    \begin{tikzpicture}
        \filldraw[black](0,0)circle (1.8pt) --(0.6,0);
        \draw[dashed](0.6,0) -- (2.6,0)node{$\circ$}  --(2.6,0);
        \draw[decorate, decoration={snake, segment length=1mm, amplitude=0.25mm}](2.6,0) -- (14.99,0);
        \draw[->](14.99,0) -- (15,0)node[right]{$ p $};
        \draw(0.6,0)node{$\circ$};
        \node[pin={[pin distance=0.2cm, pin edge={black}]{$1$}}] at (0,0) {};
        \node[pin={[pin distance=0.2cm, pin edge={black}]{$\mathfrak{p}$}}] at (0.6,0) {};
        \node[pin={[pin distance=0.2cm, pin edge={black}]{$\Bar{p}$}}] at (2.6,0) {};
    \end{tikzpicture}
    \captionsetup{singlelinecheck=off}
\caption[foo bar]{Some infimum of $\mathfrak{p}$ and $\Bar{p}$, taken from their definitions in Lemma \ref{Lemmapbarra}. From top to bottom: 
  \begin{itemize}[label={\ding{226}}]
    \item  Under assumption \ref{As2}, and  $n=1$, $m_1=\frac{81}{50}$,  $m_2=-\frac{149}{100}$, $m_3=\frac{33}{20}$, $\alpha=\frac{587}{100}$ and $\beta=\frac{63}{25}$, we have $\mathfrak{p}=\frac{n}{2}\left(m_2-m_1+\alpha\right)=\frac{69}{50}$ and $\Bar{p}=m_3(n+2)(n+1)=\frac{99}{10}$;
    \item  Under the assumption \ref{As1}, and $n=5$, $m_1=-\frac{3}{50}$, $m_2=\frac{6}{5}$, $m_3=-\frac{143}{100}$, $\alpha=\frac{163}{100}$ and $\beta=\frac{169}{100}$, we have $\mathfrak{p}=\frac{n}{2}\left(m_2-m_1+\alpha\right)=\frac{289}{40}$ and $\Bar{p}=m_2(n+2)(n+1)=\frac{252}{5}$;
    \item  Under the assumption \ref{As1}, and  $n=4$, $m_1=-\frac{187}{100}$, $m_2=-\frac{89}{100}$, $m_3=-\frac{181}{100}$, $\alpha=\frac{353}{100}$ and $\beta=\frac{19}{50}$, 
    we have $\mathfrak{p}=\frac{n}{2}\left(m_2-m_1+\alpha\right)=\frac{451}{50}$ and $\Bar{p}=n\alpha=\frac{353}{25}$;
    \item  Under the assumption  \ref{As2}, and $n=4$, $m_1=\frac{47}{50}$, $m_2=\frac{6}{25}$, $m_3=-\frac{63}{50}$, $\alpha=\frac{179}{100}$ and $\beta=\frac{119}{50}$, we have $\mathfrak{p}=\frac{n}{2}\left(m_2-m_1+\alpha\right)=\frac{109}{50}$ and $\Bar{p}=n\beta=\frac{238}{25}$.
  \end{itemize}}\label{prova5}
\end{figure}
\subsection{Improving Theorem 1.1 in \cite{wang2023blow} and addressing an open question in Remark 1.2 of \cite{wang2023blow}}\label{SubsectionImproving} Herein we want to  compare Theorem \ref{ThBlowUp} and \cite[Theorem 1.1]{wang2023blow}, both dealing with blow-up solutions to attraction-repulsion models with logistics. First, in order to have consistency between these results, we have to fix $m_2=1$ in \ref{BU}. Additionally, for the ease of the reader, let us also rephrase the related blow-up assumptions:
\begin{equation}\label{CFV}
\text{Blow-up conditions in Theorem \ref{ThBlowUp}:} \qquad 
\alpha>\beta \quad \text{and} \quad 
\begin{cases}
\alpha > \max\bigl\{m_1+\frac{2}{n}k-1, k-1\bigr\} & \text{ if } m_1 \geq 0,\\
\alpha > \max\bigl\{\frac{2}{n}k-1, k-1\bigr\} & \text{ if } m_1 < 0,\\
\end{cases}
\end{equation}
and
\begin{equation}\label{wang}
\text{Blow-up conditions in \cite[Theorem 1.1]{wang2023blow}:} \qquad 
\alpha>\beta \quad \text{and} \quad 
\begin{cases}
\alpha > \max\bigl\{m_1+\frac{2}{n}k-1, k-1\bigr\} & \text{ if } m_1 >1,\\
\alpha > \max\bigl\{\frac{2}{n}k, k-1\bigr\} & \text{ if } m_1 \leq 1.\\
\end{cases}
\end{equation}
We easily note that if $m_1 \geq 1$ the two conditions \eqref{CFV} and \eqref{wang} coincide. Now, we analyze the cases $m_1 \leq 0$ and $0<m_1<1$, separately.
\begin{itemize}[label={\ding{226}}]
\item  Case $m_1 \leq 0$.  By comparing  $\alpha > \max\bigl\{\frac{2}{n}k-1, k-1\bigr\}$ in \eqref{CFV} and $\alpha > \max\bigl\{\frac{2}{n}k, k-1\bigr\}$ in \eqref{wang}, we observe that \eqref{CFV} provides a larger range of values of $\alpha$ for which blow-up occurs than \eqref{wang} does whenever $n\in\{1,2\}$ or $n\geq 3$ provided $1<k<\frac{n}{n-2}$.
\item  Case $0<m_1<1$. The conditions above become $\alpha > \max\bigl\{m_1+\frac{2}{n}k-1, k-1\bigr\}$ and $\alpha > \max\bigl\{\frac{2}{n}k, k-1\bigr\}$, respectively. In particular, thanks to the fact that $m_1<1$, also in this situation a sharper condition is achieved for $n\in\{1,2\}$ or $n\geq 3$, under the assumptions $1<k<\frac{n}{n-2}$ and $0<m_1<\frac{k(n-2)}{n}$ or $1<k<\frac{n}{n-2}$ and $\frac{k(n-2)}{n}<m_1<1$.
\end{itemize}
From the above analysis, it is seen that Theorem \ref{ThBlowUp} improves \cite[Theorem 1.1]{wang2023blow}; additionally,  it establishes that \cite[(1.15), Theorem 1.1]{wang2023blow} is not optimal, so giving an answer to  an open question left in \cite[Remark 1.2]{wang2023blow}. 
\subsection{On the automatic applicability of Theorems \ref{maintheorem1} and \ref{maintheorem2} in some related chemotaxis contexts}\label{p0T}
We can observe what follows:
\begin{itemize}[label={\ding{226}}]
\item Once the blow-up constrains in \ref{BU} are accomplished, and taking into account $m_2=m_3>0$,  assumptions \ref{condizioneprincipale} and \ref{As1} are immediately satisfied; subsequently, Theorems \ref{maintheorem1} and \ref{maintheorem2} are applicable to unbounded solutions to model \eqref{problem3x3}.
\item Theorems \ref{maintheorem1} and \ref{maintheorem2} can also be used in models close to \eqref{problem3x3}, for which unbounded solutions can be detected (see $\S$\ref{Intro}); in particular, for attraction-repulsion (linear and nonlinear)  models with or without logistic and general production laws (\cite{wang2023blow,LiuLiBlowUpAttr-Rep}), or for only attraction ones (\cite{wang2021boundedness,black2021relaxed,WinklerNoNLinearanalysisSublinearProduction,TANAKA2021,yi2021blow}). 
\end{itemize}
\subsection{An open problem} As far as we know, establishing conditions ensuring blow-up solutions for the case $m_2\neq m_3$ in model \eqref{problem3x3}, is still an open problem; in particular, we will give some details on related technical difficulties connected to this issue in Remark \ref{Casem2neqm3}.
\section{Local existence and necessary parameters}\label{SectionLocalExistence}
By an adaption of standard reasoning in the frame of the fixed-point theorem, we can show the following result on local existence and extensibility of classical solutions to \eqref{problem3x3}.
\begin{lemma}\label{LocalExistence}
 Let $\Omega$ be a bounded and smooth domain of $\R^n$, with $n\geq 1$, $\chi,\xi,\lambda,\mu>0$, $m_1,m_2,m_3\in \R$, $k>1$, and let $f_i$ and $u_0$ comply with \eqref{f_i_u0_GENERALE}. Then there exist $\TM \in (0,\infty]$ and a unique solution $(u,v,w)$ to problem \eqref{problem3x3}, defined in $\Omega \times (0,\TM)$ and such that 
\[
u \in C^0(\bar{\Omega} \times [0,\TM)) \cap C^{2,1}(\bar{\Omega} \times (0,\TM))\;\textrm{and}\;
v, w \in \displaystyle \bigcap_{q>n} L^{\infty}_{loc}((0,\TM); W^{1,q}(\Omega)) \cap C^{2,0}(\bar{\Omega} \times (0,\TM)). 
\]
Additionally, one has $u\geq 0$ in $\Omega \times (0,\TM)$, 
\begin{equation}\label{boundednessMass}
  \int_{\Omega} u \leq M:=\max\{M_0, C\} \quad \text{for all } t \in [0,\TM),   
\end{equation} 
where $M_0=\displaystyle\int_\Omega u_0(x)\,dx$ and $C:= \left(\frac{\lambda}{\mu} |\Omega|^{k-1}\right)^{\frac{1}{k-1}}$,
and 
\begin{equation*}
\text{if} \quad \TM<\infty, \quad \text{then} \quad \limsup_{t \to \TM} \|u(\cdot,t)\|_{L^{\infty}(\Omega)}=\infty.
\end{equation*}
Furthermore, if $u_0$ satisfies also the symmetrical assumptions in \eqref{f_iDef} and $\Omega=B_R(0)$, with some $R>0$,  then $u, v$ and $w$ are radially symmetric with respect to $|x|$ in $\Omega \times (0, \TM)$.
\begin{proof}
The proof can be achieved by following well established results: for instance, we refer the interested reader to \cite{cieslak2008finite, nagai1995blow,
wang2016quasilinear, winkler2010boundedness}. In particular, an integration of the first equation in \eqref{problem3x3} and an application of the H\"{o}lder inequality give bound \eqref{boundednessMass}.
\end{proof}
\end{lemma}
Let us start with the following technical
\begin{lemma}\label{Lemmapbarra}  
Let $n\in \N$, $m_1, m_2, m_3, \alpha, \beta$ be as in Assumptions \ref{generalassumptions}. Moreover, for $\mathfrak{p}>\frac n2(m_2-m_1+\alpha)$ let 
 \begin{equation*}\label{ConstantForTechincalInequality_Barp}
 \bar{p}>\max 
\begin{Bmatrix}
 \mathfrak{p} \vspace{0.1cm}   \\ 1-m_1(n+2) \vspace{0.1cm} 
 \\ 1-m_2 \vspace{0.1cm}
  \\ 1-m_3 \vspace{0.1cm}
  \\ 2-m_1-\frac{2}{n} \vspace{0.1cm} 
 \\n \alpha \vspace{0.1cm}
  \\n \beta \vspace{0.1cm}\\ 
  m_2(n+2)(n+1) \vspace{0.1cm}  \\ m_3(n+2)(n+1)
 \end{Bmatrix},
 \end{equation*}
\begin{equation*}
    \sigma\coloneqq\frac{2(p+m_2+\alpha-1)}{p+m_1-1}, \quad \hat{\sigma}\coloneqq\frac{2p}{p+m_1-1}, \quad \gamma\coloneqq\frac{\frac{m_1-m_2-\alpha}{2}+\frac pn +\frac{m_2+\alpha-1}{n}}{\frac{m_1-m_2-\alpha}{2}+\frac pn}, \quad \delta \coloneqq\frac{p+m_2+\alpha-1}{p}.
\end{equation*}
Then for all $p \geq  \bar{p}$ these relations hold

 \begin{minipage}{0.6\textwidth}
 \begin{subequations}
 \begin{tabularx}{\textwidth}{XX}
 \begin{equation}\label{2}
     p>\frac n2 (1-m_1),
 \end{equation}
 &
 \begin{equation}\label{theta}
     0<\theta\coloneqq\frac{\frac{p+m_1-1}{2\mathfrak{p}}-\frac{p+m_1-1}{2(p+m_2+\alpha-1)}}{\frac{p+m_1-1}{2\mathfrak{p}}+\frac 1n -\frac 12}<1,
 \end{equation}
 \\
 \begin{equation}\label{y1}
     0<\frac{\sigma\theta}{2}<1,
 \end{equation}
&
 \begin{equation}\label{thetacap}
     0<\hat{\theta}=\frac{\frac{p+m_1-1}{2}-\frac{p+m_1-1}{2p}}{\frac{p+m_1-1}{2}+\frac 1n -\frac 12}<1,
 \end{equation}
  \\
 \begin{equation}\label{hatsigmatheta2}
     0<\frac{\hat{\sigma}\hat{\theta}}{2}<1,
 \end{equation}
 &
 \begin{equation}\label{esp1}
     0<\frac{p+m_3-1}{p+m_2+\alpha-1}<1,
 \end{equation}
  \\
 \begin{equation}\label{esp2}
     0<\frac{\beta}{p+m_2+\alpha-1}<1,
 \end{equation}
 &
 \begin{equation}\label{sommaesp}
     \frac{p+m_3-1}{p+m_2+\alpha-1}+\frac{\beta}{p+m_2+\alpha-1}<1,
 \end{equation}
 \\
 \begin{equation}\label{esp3}
     0<\frac{p}{p+m_2+\alpha-1}<1,
 \end{equation}
 &
 \begin{equation}\label{gamma>delta>1}
     \gamma > \delta > 1,
 \end{equation}
 \\
 \begin{equation}\label{thetabar}
    0<\Bar{\theta}\coloneqq\frac{\frac{p+m_1-1}{2p}-\frac{p+m_1-1}{2(p+m_2+\alpha-1)}}{\frac{p+m_1-1}{2p}+\frac 1n -\frac 12}<1,
 \end{equation}
 &
 \begin{equation}\label{barsigmatheta2}
     0<\frac{\sigma\Bar{\theta}}{2}<1,
 \end{equation}
 \begin{equation}\label{Y3}
     0<\frac{p+m_3-1}{p}<1.
 \end{equation}
 \end{tabularx}
 \end{subequations}
 \end{minipage}
\begin{proof}
To show our relations, we will need $p>1-m_1$ and $p>1+\beta-m_2-\alpha$. As to the first inequality, due to \ref{condizioneprincipale} and the restriction on $\mathfrak{p}$, we have that $p \geq \bar{p}>\mathfrak{p}>1\geq 1-m_1$ for $m_1\geq0$, whilst for $m_1<0$ it suddenly derives from $p\geq \bar{p}>1-m_1(n+2)$. For the second lower bound,  assumptions \ref{As1} or \ref{As2}, together with the definition of $\bar{p}$, give $p\geq \bar{p}>1-m_3>1+\beta-m_2-\alpha$. Besides, $m_2+\alpha>1$ is automatically true through \ref{As1}, or alternatively by means of assumption \ref{condizioneprincipale} in conjunction with one among \ref{As2}, \ref{As3} or \ref{As4}. The same condition $m_2+\alpha>1$, $p\geq \bar{p}$ and the restriction on $\mathfrak{p}$ and $p>1-m_1$ ensure relations \eqref{2}, \eqref{theta}, \eqref{y1}, \eqref{thetabar} and \eqref{barsigmatheta2}.
Moreover, from \eqref{2}, $p>1-m_1$ and $p>2-m_1-\frac{2}{n}$, we obtain \eqref{thetacap}
and by using also $m_1>1-\frac 2n$, we get \eqref{hatsigmatheta2}.  On the other hand, restrictions \eqref{esp1}, \eqref{esp2}, \eqref{sommaesp}, \eqref{esp3} and \eqref{Y3} come from the definition of $\bar{p}$ in conjunction with $m_2+\alpha>m_3$, $p>1+\beta-m_2-\alpha$, $m_2+\alpha>m_3+\beta$ and $m_2+\alpha>1$, respectively.
Finally, from \ref{condizioneprincipale} it follows that $m_2+\alpha>m_1$, which combined with $m_2+\alpha>1$ gives \eqref{gamma>delta>1}. 
\end{proof}
\end{lemma}
\begin{remark}\label{RemarkCasoLimite}
 For reasons which will be exploited later on, and precisely in Lemma \ref{lemmateo1}, it appears important to point out that assumptions \ref{As2}, \ref{As3} and \ref{As4} imply that the ratios  
 \eqref{esp1}, \eqref{sommaesp} and \eqref{Y3} in Lemma \ref{Lemmapbarra} can be also taken equal to 1.  
 \end{remark}
\section{A priori estimates and proof of Theorems \ref{maintheorem1} and \ref{maintheorem2}}
In this section we will use \ref{condizioneprincipale}--- \ref{As4}. Moreover, without explicitly computing their values, we underline that the constants $c_i$ appearing below and throughout the paper depend inter alia on $p$, are positive and their subscripts $i$ start anew in each new proof. 
\resetconstants{c}
\begin{lemma}\label{lemmateo1}
Under the hypotheses of Lemma \ref{Lemmapbarra}, let $p=\Bar{p}$ and $\mathfrak{p}$ be any of the constants therein defined. If $(u,v,w)$ is a classical solution to problem \eqref{problem3x3},  $u\in L^{\infty}((0,\TM);L^{\mathfrak{p}}(\Omega))$ and $\varphi(t)$ is the energy function
\begin{equation*}
    \varphi(t)\coloneqq\frac 1p \into (u+1)^p \quad \text{on $(0,\TM)$},
\end{equation*}
then there exist $\const{GiuA},\const{GiuB}$ such that
\begin{equation}\label{claimlemmateo1}
    \varphi'(t)\leq -\const{GiuA} \into \abs*{\nabla (u+1)^{\frac{p+m_1-1}{2}}}^2 +\const{GiuB} \quad \text{for all $t\in(0,\TM)$}.
\end{equation}
\begin{proof}
Let us differentiate the functional $\varphi(t)=\frac 1p \into (u+1)^p$. Using the first equation of \eqref{problem3x3} and the divergence theorem we have for every $t\in(0,\TM)$
\begin{equation*}
\begin{split}
    \varphi'(t)=\into (u+1)^{p-1}u_t=&\into (u+1)^{p-1}\nabla\cdot ((u+1)^{m_1-1}\nabla u)-\chi\into (u+1)^{p-1}\nabla\cdot (u(u+1)^{m_2-1}\nabla v)\\
    &+\xi\into(u+1)^{p-1}\nabla\cdot (u(u+1)^{m_3-1}\nabla w) +\lambda\into (u+1)^{p-1}u-\mu\into (u+1)^{p-1}u^k\\
    =&-(p-1)\into(u+1)^{p+m_1-3}\abs{\nabla u}^2+(p-1)\chi\into u(u+1)^{p+m_2-3}\nabla u\cdot\nabla v\\
    &-(p-1)\xi\into u(u+1)^{p+m_3-3}\nabla u\cdot \nabla w +\lambda\into (u+1)^{p-1}u-\mu\into (u+1)^{p-1}u^k.
\end{split}
\end{equation*}
For $j\in\{m_2,m_3\}$, we now define 
\begin{equation*}
    F_j(u)=\int_0^u \hat{u}(\hat{u}+1)^{p+j-3}d\hat{u},
\end{equation*}
so observing that
\begin{equation}\label{disF}
    0\leq F_j(u)\leq \frac{1}{p+j-1}\left[(u+1)^{p+j-1}-1\right].
\end{equation}
By considering the definition of $F_j(u)$ above, again the divergence theorem, the second and third equation of \eqref{problem3x3}, we have for every $t\in(0,\TM)$,
\begin{equation}\label{disprincipale}
\begin{split}
    \varphi'(t)\leq &-\const{a}\into \abs*{\nabla(u+1)^{\frac{p+m_1-1}{2}}}^2+\const{b}\into \nabla F_{m_2}(u)\cdot \nabla v -\const{c}\into \nabla F_{m_3}(u)\cdot \nabla w +\const{d}\into(u+1)^p-\const{e}\into(u+1)^{p-1}u^k\\
    =&-\const{a}\into \abs*{\nabla(u+1)^{\frac{p+m_1-1}{2}}}^2-\const{b}\into F_{m_2}(u)\Delta v +\const{c}\into F_{m_3}(u)\Delta w +\const{d}\into(u+1)^p-\const{e}\into(u+1)^{p-1}u^k.
\end{split}
\end{equation}
Let us now specify how each of the constrains in Assumptions \ref{generalassumptions} takes part in our computation. 

First, from \ref{condizioneprincipale} we have $\mathfrak{p}>1$; this makes meaningful our assumption  $u\in L^{\infty}((0,\TM);L^{\mathfrak{p}}(\Omega))$. (Recall that  $u\in L^{\infty}((0,\TM);L^1(\Omega))$ is always met by \eqref{boundednessMass}.)

Now, by exploiting \ref{As1} and \eqref{disF}, we can see that from \eqref{disprincipale}, if we neglect the nonpositive terms we get on $(0,\TM)$
\begin{equation}\label{phi1}
    \varphi'(t)\leq -\const{a}\into \abs*{\nabla(u+1)^{\frac{p+m_1-1}{2}}}^2 +\const{f}\into (u+1)^{p+m_2+\alpha-1}+\const{g}\into (u+1)^{p+m_3-1}\into (u+1)^\beta+\const{d}\into (u+1)^p.
\end{equation}
As to the third term, by using twice H\"{o}lder's inequality (recall \eqref{esp1} and \eqref{esp2}), we obtain, for every $t\in(0,\TM)$,
\begin{equation}\label{InGN}
    \const{g}\into (u+1)^{p+m_3-1}\into (u+1)^\beta \leq \const{i}\tonda*{\into (u+1)^{p+m_2+\alpha-1}}^{\frac{p+m_3-1}{p+m_2+\alpha-1}}\tonda*{\into (u+1)^{p+m_2+\alpha-1}}^{\frac{\beta}{p+m_2+\alpha-1}}.
\end{equation}
Moreover, since for any $\varepsilon>0$ there is $d(\varepsilon)>0$ such that this inequality (see \cite[Lemma 4.3]{frassuviglialoro1}) 
\begin{equation*}
    A^{d_1}B^{d_2}\leq \varepsilon(A+B)+d(\varepsilon), \quad A,B\geq0, \, d_1,d_2>0, \, d_1+d_2<1,
\end{equation*}
is true, by virtue of \eqref{sommaesp}, we have that
\begin{equation}\label{termine1}
    \const{g}\into (u+1)^{p+m_3-1}\into (u+1)^\beta \leq \into (u+1)^{p+m_2+\alpha-1}+\const{l}
    \quad  \textrm{for all $t\in(0,\TM)$}.
\end{equation}
Through \eqref{esp3}, an application of Young's inequality provides 
\begin{equation}\label{termine2}
    \const{d}\into (u+1)^p\leq \const{d5}\tonda*{\into (u+1)^{p+m_2+\alpha-1}}^{\frac{p}{p+m_2+\alpha-1}}\leq \into (u+1)^{p+m_2+\alpha-1}+\const{n} \quad \textrm{on } (0,\TM).
\end{equation}
To control the term $\into (u+1)^{p+m_2+\alpha-1}$, we invoke the Gagliardo--Nirenberg and Young's inequalities, so to bound the mentioned integral with $\into |\nabla(u+1)^{\frac{p+m_1-1}{2}}|^2$. More exactly by relying on relations \eqref{theta} and \eqref{y1}, boundedness of $\into (u+1)^{\mathfrak{p}}$ and of $\into u$ provide for any $L_1>0$ 
\begin{equation}\label{secondotermine}
\begin{split}
    L_1\into (u+1)^{p+m_2+\alpha-1}=& L_1\norm*{(u+1)^\frac{p+m_1-1}{2}}^{\frac{2(p+m_2+\alpha-1)}{p+m_1-1}}_{L^{\frac{2(p+m_2+\alpha-1)}{p+m_1-1}}(\Omega)}\\ 
    \leq & \const{q}\tonda*{\norm*{\nabla(u+1)^{\frac{p+m_1-1}{2}}}^{\sigma\theta}_{L^2(\Omega)}\norm*{(u+1)^{\frac{p+m_1-1}{2}}}^{\sigma(1-\theta)}_{L^{\frac{2\mathfrak{p}}{p+m_1-1}}(\Omega)}+\norm*{(u+1)^{\frac{p+m_1-1}{2}}}_{L^{\frac{2}{p+m_1-1}}(\Omega)}^{\sigma}}\\
    \leq & \const{w}\tonda*{\into \abs*{\nabla(u+1)^{\frac{p+m_1-1}{2}}}^2}^{\frac{\sigma\theta}{2}}+\const{r} \leq  \frac{c_3}{2} \into \abs*{\nabla(u+1)^{\frac{p+m_1-1}{2}}}^2 +\const{t}  \quad  \textrm{for all } t\in(0,\TM).
\end{split}
\end{equation} 
(We underline that herein we have used the elementary inequality 
\begin{equation}\label{algebricIne}
C_1(\tau)\left(A^\tau+B^\tau\right)\leq (A+B)^{\tau} \leq C_2(\tau) (A^{\tau}+B^{\tau})\quad \textrm{for all } A,B\geq 0, \tau>0 \textrm{ and proper } C_1,C_2>0,
\end{equation}
which might tacitly be used in the next lines.)
Putting together \eqref{phi1}, \eqref{termine1}, \eqref{termine2}, \eqref{secondotermine}, we have the claim.

If assumption \ref{As2} is complied, bound \eqref{termine2} can be replaced for any $L_2>0$ with 
\begin{equation}\label{gnconas2}
\begin{split}
    L_2\into (u+1)^p = & L_2\norm*{(u+1)^{\frac{p+m_1-1}{2}}}^{\frac{2p}{p+m_1-1}}_{L^{\frac{2p}{p+m_1-1}}(\Omega)}\\
    \leq & \const{a1}\tonda*{\norm*{\nabla(u+1)^{\frac{p+m_1-1}{2}}}^{\hat{\sigma}\hat{\theta}}_{L^2(\Omega)}\norm*{(u+1)^{\frac{p+m_1-1}{2}}}^{\hat{\sigma}(1-\hat{\theta)}}_{L^{\frac{2}{p+m_1-1}}(\Omega)}+\norm*{(u+1)^{\frac{p+m_1-1}{2}}}_{L^{\frac{2}{p+m_1-1}}(\Omega)}^{\hat{\sigma}}}\\
    \leq & \const{a1}\tonda*{\into \abs*{\nabla(u+1)^{\frac{p+m_1-1}{2}}}^2}^{\frac{\hat{\sigma}\hat{\theta}}2}+\const{c1}
    \leq  \frac{c_3}{4}\into \abs*{\nabla(u+1)^{\frac{p+m_1-1}{2}}}^2+\const{b1}\quad \textrm{on } (0,\TM),
\end{split}
\end{equation}
where in this last step we used Young's inequality in conjunction with \eqref{thetacap} and \eqref{hatsigmatheta2}. 
Moreover, by taking $m_2+\alpha=m_3+\beta$ in \ref{As2} (recall Remark \ref{RemarkCasoLimite}), in estimate \eqref{InGN} the powers at the rhs have 1 as sum; subsequently, up to constants, it becomes \eqref{termine1}.
So, by considering \eqref{phi1}, \eqref{termine1}, \eqref{secondotermine} and \eqref{gnconas2}  we have the claim when either \ref{condizioneprincipale} and \ref{As1} or \ref{condizioneprincipale} and \ref{As2} are used. 

Let us show how achieving the same conclusion by applying either \ref{condizioneprincipale} and \ref{As3} or \ref{condizioneprincipale} and \ref{As4}. The main idea is 
alternatively treating in relation \eqref{phi1} the term
$$
\into (u+1)^{p+m_3-1}\into (u+1)^\beta  \quad \textrm{on } (0,\TM).
$$
More specifically, if $\beta\in \left(0,1\right]$ and we take into account the boundedness of $\int_\Omega u$ on $(0,\TM)$, we have that 
\begin{equation}\label{betaMnoUnoFrassu}
\into (u+1)^{p+m_3-1}\into (u+1)^\beta \leq \const{Giu1} \int_\Omega (u+1)^{p+m_3-1} \quad \textrm{for all } t \in (0,\TM).
\end{equation}
In turn, by relying respectively on assumption \ref{As3} or \ref{As4} (both with strict inequality), by means of Young's inequality, thanks to \eqref{Y3} and \eqref{esp1}, it is also possible to see that
\begin{equation}\label{betaMnoUnoFrassuBIS}
 \const{Giu1} \int_\Omega (u+1)^{p+m_3-1} \leq \int_\Omega (u+1)^p+\const{Giu2}\quad 
 \textrm{on } (0,\TM),
\end{equation}
or
\begin{equation}\label{betaMnoUnoFrassuTRIS}
\const{Giu1} \int_\Omega (u+1)^{p+m_3-1} \leq 
 \int_\Omega (u+1)^{p+m_2-1+\alpha}+\const{Giu3}\quad \textrm{for all } t \in (0,\TM).
\end{equation}
As before, bounds  \eqref{phi1}, \eqref{betaMnoUnoFrassu}, \eqref{betaMnoUnoFrassuBIS}, or alternatively \eqref{betaMnoUnoFrassuTRIS},   \eqref{gnconas2} and \eqref{secondotermine}, lead to the same conclusion. 
(For the limit cases in \ref{As3} or \ref{As4}, namely $m_3=1$ or $m_2+\alpha=m_3$, by relying again on Remark \ref{RemarkCasoLimite}, we can directly exploit estimate \eqref{betaMnoUnoFrassu}, without using Young's inequality, and conclude as above.)
\end{proof}
\end{lemma}
The next step consists in ensuring some time independent estimate of $u$ in the $L^{\bar{p}}(\Omega)$-norm.
\resetconstants{c}
\begin{lemma}\label{lemmateo12}
Let Lemma \ref{lemmateo1} be true. Then $u\in L^{\infty}((0,\TM);L^{\Bar{p}}(\Omega))$. 
\begin{proof}
With a view to inequality \eqref{claimlemmateo1}, as already done in \eqref{gnconas2}, a further application of the Gagliardo--Nirenberg inequality, supported by \eqref{thetacap}, leads also thanks to \eqref{algebricIne} to 
\begin{equation*}
\begin{split}
    \into (u+1)^p 
    \leq \const{d1}\tonda*{\norm*{\nabla(u+1)^{\frac{p+m_1-1}{2}}}^{\hat{\sigma}\hat{\theta}}_{L^2(\Omega)}+1}
    \leq  \const{e1}\tonda*{\into \abs*{\nabla(u+1)^{\frac{p+m_1-1}{2}}}^{2}+1}^{\frac{{\hat{\sigma}\hat{\theta}}}{2}} \quad \textrm{on } (0,\TM),
\end{split}
\end{equation*}
or, equivalently
\begin{equation}\label{dis3}
    -\const{f1}\into \abs*{\nabla(u+1)^{\frac{p+m_1-1}{2}}}^{2}\leq -\tonda*{\into (u+1)^p}^{\frac{2}{\hat{\sigma}\hat{\theta}}}+\const{f1} \quad \text{for every $t\in(0,\TM)$}.
\end{equation}
Finally, by using relations \eqref{claimlemmateo1} and \eqref{dis3}, we arrive at this initial problem
\begin{equation*}
\begin{dcases}
    \varphi'(t)\leq \const{k1}-\const{l1}\varphi(t)^{\frac{2}{\hat{\sigma}\hat{\theta}}} \quad \textrm{for every } t\in(0,\TM),\\
    \varphi(0)=\frac{1}{p}\into (u_0+1)^p,
\end{dcases}
\end{equation*}
which ensures that $\into u^{\bar{p}}\into u^p\leq \into (u+1)^p\leq p \max\left\{\varphi(0),\tonda*{\frac{\const{k1}}{\const{l1}}}^\frac{\hat{\sigma}\hat{\theta}}{2}\right\}$ for all $t<\TM$. 
\end{proof}
\end{lemma}
By taking advantage from the previous lemma, let us show the uniform-in-time boundedness of $u$.
\resetconstants{c}
\begin{lemma}\label{lemmateo13}
Under the hypotheses of Lemma \ref{lemmateo12}, $u\in L^{\infty}((0,\TM);L^{\infty}(\Omega))$.
\begin{proof}
With the same nomenclature used by Tao and Winkler, $u$ also classically solves in $\Omega \times (0,T_{max})$ problem (A.1) of \cite[Appendix A]{TaoWinkParaPara} for 
\begin{equation*}
D(x,t,u)=(u+1)^{m_1-1},\quad f(x,t)=\chi u(u+1)^{m_2-1}\nabla v-\xi  u(u+1)^{m_3-1}\nabla w,\quad g(x,t)=\lambda u-\mu u^{k}. 
\end{equation*}
In particular, also taking into account the boundary condition on $v$ and $w$, we can see that (A.2)--(A.5) are met. On the other hand, for any $\lambda,\mu>0$ and $k>1$, it holds that $\lambda u-\mu u^{k}$ has a positive maximum $L$ at $u_M=\left(\frac{\lambda}{k \mu}\right)^{\frac{1}{k-1}}$, so that from $g(x,t)\leq L$ in $\Omega \times (0,T_{max})$ the second inclusion of (A.6) is accomplished for any choice of $q_2$.  As to (A.7)--(A.10), let us first define the quantities 
\begin{equation*}
    l_1(q_1)=1-m_1\frac{(n+1)q_1-(n+2)}{q_1-(n+2)},\quad l_2(q_2)=1-m_1\frac{1}{1-\frac{n q_2}{(n+2)(q_2-1)}}\quad \textrm{and}\quad l_3=\frac{n}{2}\left(1-m_1\right). 
\end{equation*}
Recalling the definition and properties of $p=\bar p$, we have $p>1-m_1(n+2)$ and henceforth  for any $m_1\in \R$, it holds  $1-m_1(n+2)\geq l_1((n+2)(n+1))$ and  there exists $\tilde{q}_2$ large enough so to have $1-m_1(n+2)\geq l_2(\tilde{q}_2)$. Subsequently, \cite[(A.7)]{TaoWinkParaPara} (i.e. $u\in L^\infty((0,T_{max});L^p(\Omega))$, that is $\into u^p \leq \const{Miso}$) is obviously accomplished,  \cite[(A.8), Lemma A.1.]{TaoWinkParaPara}  is fulfilled for $q_1=(n+2)(n+1)$, \cite[(A.9), Lemma A.1.]{TaoWinkParaPara}  for sufficiently large $q_2$, whereas  \cite[(A.10), Lemma A.1.]{TaoWinkParaPara} is directly true thanks to \eqref{2}. Let us dedicate to the first inclusion  of \cite[(A.6)]{TaoWinkParaPara}. Starting from the gained bound of $u$, let us exploit elliptic regularity results applied to the second and third equation of system \eqref{problem3x3}; in particular, we only analyze $-\Delta v=f_1(u)-\frac{1}{|\Omega|}\int_\Omega f_1(u)$, being the case for $w$ equivalent. Let us observe that from relation \eqref{CondizioniSuFAlessandro} on $f_1$, and $p=\bar{p}>n\alpha$ naturally implying $(1+u)^\alpha\leq (1+u)^p$, we have for all  $t\in (0,\TM)$
\begin{equation*}
\int_\Omega |f_1(u)-\frac{1}{|\Omega|}\int_\Omega f_1(u)|^{\frac{p}{\alpha}}\leq \const{G1}\int_\Omega (1+u)^{p}+\const{G2}\int_\Omega \left(\int_\Omega (1+u)^{\alpha }\right)^{\frac{p}{\alpha}}\leq \const{G1}\int_\Omega (1+u)^{p}+\const{G2}\int_\Omega \left(\int_\Omega (1+u)^{p }\right)^{\frac{p}{\alpha}}\leq \const{G10},
\end{equation*}
which gives $f_1(u)-\frac{1}{|\Omega|}\int_\Omega f_1(u)\in L^\infty((0,T_{max});L^{\frac{p}{\alpha}}(\Omega))$, in turn $v\in L^\infty((0,T_{max});W^{2,\frac{p}{\alpha}}(\Omega))$ and, hence, through the Sobolev embeddings $\nabla v \in L^\infty((0,T_{max});W^{1,\frac{p}{\alpha}}(\Omega)) \xhookrightarrow{}  L^\infty((0,T_{max});L^{\infty}(\Omega))$. Consequently, thanks to the H\"{o}lder inequality (recall that from  Lemma \ref{Lemmapbarra} one has that  $p>m_2(n+2)(n+1)$), by using the uniform-in-time boundedness of $u$ in $L^{\Bar{p}}(\Omega)$ we can write on $(0,T_{max})$
\begin{equation*}
\begin{split}
\int_\Omega |u(u+1)^{m_2-1}\nabla v|^{(n+2)(n+1)}&
\leq  \lVert \nabla v (\cdot, t) \rVert_{L^\infty(\Omega)}^{(n+2)(n+1)}|\Omega|^\frac{p-m_2(n+2)(n+1)}{p} \left(\int_\Omega (u+1)^{p}\right)^\frac{m_2(n+2)(n+1)}{p}\leq \const{G6}.
\end{split}
\end{equation*}
Reasoning in a similar way on the third equation of \eqref{problem3x3}, we have $\nabla w\in L^\infty((0,T_{max});L^{\infty}(\Omega))$,
\begin{equation*}
\begin{split}
\int_\Omega |u(u+1)^{m_3-1}\nabla w|^{(n+2)(n+1)}\leq \const{G7} \quad \textrm{for all } t \in (0,\TM),
\end{split}
\end{equation*}
and as a consequence
\begin{equation*}
f=\chi u(u+1)^{m_2-1}\nabla v -\xi u(u+1)^{m_3-1}\nabla w \in L^\infty((0,T_{max});L^{(n+2)(n+1)}(\Omega)). 
\end{equation*}
Since all the hypotheses of \cite[Lemma A.1.]{TaoWinkParaPara} are fulfilled, we have the claim. 
\end{proof} 
\end{lemma}
Now we are in a position to prove our first results.
\resetconstants{c}
\subsubsection*{Proof of Theorem \ref{maintheorem1}:}
Let $(u,v,w)$ be a given blow-up solution at some finite time $\TM$ to problem \eqref{problem3x3}. If $u$ was not unbounded in some $L^{\mathfrak{p}}(\Omega)$-norm, Lemma  \ref{lemmateo13} would imply the uniform-in-time boundedness of $u$, contradicting hypothesis \eqref{blowup}.
\qed
\subsubsection*{Proof of Theorem \ref{maintheorem2}:} By making use of \eqref{termine1} only or altogether \eqref{betaMnoUnoFrassu}, \eqref{betaMnoUnoFrassuBIS} and \eqref{betaMnoUnoFrassuTRIS}, we observe that bound \eqref{phi1} can essentially be reorganized as
\begin{equation}\label{phi2}
    \varphi'(t)\leq -\const{a}\into \abs*{\nabla(u+1)^{\frac{p+m_1-1}{2}}}^2 +\const{fK}\into (u+1)^{p+m_2+\alpha-1}+
    \const{hK}\into (u+1)^p+\const{hK2} \quad \textrm{on } (0,\TM).
\end{equation}
From the one hand, as already done in \eqref{termine2}, up to a constant the Young inequality yields
\begin{equation}\label{BlowUpTime1}
\into (u+1)^p\leq \int_\Omega (u+1)^{p+m_2+\alpha-1}+\const{zG} \quad \textrm{for all } t \in (0,\TM);
\end{equation}
moreover, by recalling \eqref{thetabar} and \eqref{barsigmatheta2}, we can derive again through the Gagliardo--Nirenberg and Young's inequalities, this bound valid for any $K>0$
\begin{equation}\label{phi1teo2}
\begin{split}
    K\into (u+1)^{p+m_2+\alpha-1} = & K \norm*{(u+1)^\frac{p+m_1-1}{2}}^{\frac{2(p+m_2+\alpha-1)}{p+m_1-1}}_{L^{\frac{2(p+m_2+\alpha-1)}{p+m_1-1}}(\Omega)}\\ 
    \leq & \const{q2}\tonda*{\norm*{\nabla(u+1)^{\frac{p+m_1-1}{2}}}^{\sigma\Bar{\theta}}_{L^2(\Omega)}\norm*{(u+1)^{\frac{p+m_1-1}{2}}}^{\sigma(1-\Bar{\theta})}_{L^{\frac{2p}{p+m_1-1}}(\Omega)}+\norm*{(u+1)^{\frac{p+m_1-1}{2}}}_{L^{\frac{2p}{p+m_1-1}}(\Omega)}^{\sigma}}\\
    \leq & \const{w2}\tonda*{\into \abs*{\nabla(u+1)^{\frac{p+m_1-1}{2}}}^2}^{\frac{\sigma\Bar{\theta}}{2}}\tonda*{\into (u+1)^p}^{\frac{p+m_1-1}{2p}\sigma(1-\Bar{\theta})}+\const{5r}\tonda*{\into (u+1)^p}^{\frac{p+m_1-1}{2p}\sigma} \\
    \leq &  \const{a} \into \abs*{\nabla(u+1)^{\frac{p+m_1-1}{2}}}^2 +\const{t2}\tonda*{\into (u+1)^p}^{\gamma}+\const{s2}\tonda*{\into (u+1)^p}^{\delta} \quad \textrm{on }  (0,\TM).
\end{split}
\end{equation}
Finally, by inserting relations \eqref{BlowUpTime1} and \eqref{phi1teo2} into \eqref{phi2}, and by using \eqref{phi1teo2} with $K= \const{fK} + \const{hK}$, we obtain for proper positive $\mathcal{A}, \mathcal{B}$ and $\mathcal{C}$
\begin{equation}\label{phi3}
    \varphi'(t)\leq \mathcal{A}\varphi^\gamma(t) +\mathcal{B}\varphi^\delta(t) + \mathcal{C} \quad \textrm{on }  (0,\TM).
\end{equation}
 Since $p=\bar{p}>\mathfrak{p}$, from Theorem \ref{maintheorem1} we know that $\limsup_{t\rightarrow T_{max}}\frac{1}{p}\int_\Omega (u+1)^p=\infty.$ On the other hand, since $\varphi(t)$ satisfies relation  \eqref{phi3} for  any $0<t<T_{max}$, the function $\Psi(\xi)=\mathcal{A}\xi^\gamma+\mathcal{B}\xi^\delta+\mathcal{C}$ obeys  the Osgood criterion \eqref{OsgoodCriterion}. Thereafter, by integrating \eqref{phi3} between $0$ and $T_{max}$,  we obtain estimate \eqref{stimatempobu}, and the first conclusion is achieved. 
 
As to the derivation of the explicit expression for the lower bound $T$, let us reduce \eqref{phi3} as follows: from $|\Omega|\leq \int_\Omega (u+1)^p=p \varphi(t)$, we can estimate $\mathcal{C}$ in relation \eqref{phi3} as
\begin{equation*} 
\mathcal{C} \leq p\frac{\mathcal{C}}{|\Omega|}\varphi(t)=:\bar{\mathcal{C}}\varphi(t),
\end{equation*} 
so that \eqref{phi3}  can be rewritten in this form:
\begin{equation}\label{ldiff}
 \varphi'(t)\leq \mathcal{A}\varphi^\gamma(t) +\mathcal{B}\varphi^\delta(t) + \bar{\mathcal{C}}\varphi(t)\quad \textrm{on }  (0,\TM).
\end{equation} 
Now, since $\varphi$ blows up at finite time $T_{max}$, there exists a time $t_1 \in [0, T_{max})$ such that 
\begin{equation*}
\varphi(t) \geq \varphi(0) \quad  \textrm{for all } \, t \geq t_1 \in[0,T_{max}).
\end{equation*}
From  $\gamma >\delta >1 $ (recall \eqref{gamma>delta>1}), we can estimate the second and third  term on the rhs of \eqref{ldiff} by means of $\varphi^{\gamma}$:
\begin{equation} \label{Phi^delta < Phi^gamma}
\varphi^{\delta}(t) \leq  \varphi(0)^{\delta - \gamma} \varphi^{\gamma}(t)\quad \textrm{ and }
\quad \varphi(t) \leq \varphi(0)^{1- \gamma} \varphi^{\gamma}(t)\quad \textrm{for all } \, t \geq t_1 \in[0,T_{max}).
\end{equation}
By plugging expressions \eqref{Phi^delta < Phi^gamma} into \eqref{ldiff} we obtain for
\begin{equation*}
\mathcal{D}=\mathcal{A}+ \mathcal{B} \varphi(0)^{\delta -\gamma} + \bar{\mathcal{C}} \varphi(0)^{1 - \gamma},
\end{equation*}
\begin{equation}\label{ldiff bis}
\varphi'(t)\leq \mathcal{D} \varphi^{\gamma}(t) \quad \textrm{for all } \, t \geq t_1 \in[0,T_{max}), 
\end{equation}
so that an integration of \eqref{ldiff bis} on $(t_1,T_{max})$ yields this explicit lower bound for $T_{max}$:
\begin{align*}
T = \frac{\varphi(0)^{1-\gamma}}{\mathcal{D} (\gamma -1)}= \int_{\varphi(0)}^{\infty} \frac{d\tau} {\mathcal{D}\tau^{\gamma}}\leq \int_{t_1}^{T_{max}}  d\tau\leq \int_0^{T_{max}} d\tau =  T_{max}. 
\end{align*}
\qed
\begin{remark}
    For completeness, we observe that it is also possible to avoid estimate \eqref{BlowUpTime1}; indeed, in 
    relation \eqref{phi2} the term $\int_\Omega (u+1)^p$ is directly $p\varphi(t)$, so that \eqref{phi3} would read
\begin{equation*}
    \varphi'(t)\leq \tilde{\mathcal{A}}\varphi^\gamma(t) +\tilde{\mathcal{B}}\varphi^\delta(t) + \tilde{\mathcal{C}}\varphi(t) +\tilde{\mathcal{D}}\quad \textrm{on }  (0,\TM).
\end{equation*}
    
    \end{remark}
\resetconstants{c}
\section{Finite-time blow-up to a simplified version of problem \eqref{problem3x3}}\label{BlowUp}
This section is dedicated to prove finite time blow-up for solutions to problem \eqref{problem3x3} in a more specific case; in our computations we will be inspired by \cite{TANAKA2021,wang2023blow}, where respectively blow-up is established in a model with only attraction, nonlinear diffusion and sensitivity and  logistic term, and in an attraction-repulsion one, with nonlinear diffusion but linear sensitivities and logistics.
\subsection{Detecting unbounded solutions to problem \eqref{problem3x3} for $m_2=m_3>0$}
Let us fix $m_2=m_3>0$ in model \eqref{problem3x3}, and in turn let us set $z=\chi v -\xi w$, $m(t)=\chi m_1(t)-\xi m_2(t)$ and $f(u)=\chi f_1(u)-\xi f_2(u)$, being $m_1(t)$ and $m_2(t)$ defined in \eqref{Definitionm1andme}; in this way, problem \eqref{problem3x3} itself is reduced into 
\begin{equation}\label{ReducedProblem}
\begin{cases}
u_t= \nabla \cdot ((u+1)^{m_1-1}\nabla u - u(u+1)^{m_2-1}\nabla z)+ \lambda u -\mu u^k  & \text{ in } \Omega \times (0,T_{max}),\\
0=\Delta z-m(t)+f(u)  & \text{ in } \Omega \times (0,T_{max}),\\
u_{\nu}=z_{\nu}=0 & \text{ on } \partial \Omega \times (0,T_{max}),\\
u(x,0)=u_0(x) & x \in \bar\Omega,\\
\int_\Omega z(x,t)dx = 0 & \textrm{for all } t \in (0,\TM).
\end{cases}
\end{equation}
In particular, if we confine our study to radially symmetric cases, by setting $r:=|x|$ and by considering $\Omega=B_R(0)\subset \R^n$, $n\geq 1$ and some $R>0$,  the radially symmetric local solution 
\[
(u,z)=(u(r,t), z(r,t))
\]
to model \eqref{ReducedProblem}  solves the following scalar problem
\begin{equation}\label{RadialProblem}
\begin{cases}
r^{n-1} u_t= (r^{n-1} (u+1)^{m_1-1}u_r)_r - (r^{n-1} u(u+1)^{m_2-1} z_r)_r  + \lambda 
r^{n-1} u -\mu r^{n-1} u^k  & r \in (0,R), t \in (0,T_{max}),\\
0= (r^{n-1} z_r)_r - r^{n-1} m(t) + r^{n-1} f(u) & r \in (0,R), t \in (0,T_{max}),\\
u_r=z_r=0 & r=R, t \in (0,T_{max}),\\
u(r,0)=u_0(r) & r \in (0,R),\\
\int_0^R r^{n-1} z(r,t) \,dr = 0 & \textrm{for all } t \in (0,\TM).
\end{cases}
\end{equation}
In the same spirit of \cite{JaLu}, we introduce the mass accumulation function
\begin{equation}\label{W}
U(s,t):= \int_0^{s^{\frac{1}{n}}} \rho^{n-1} u(\rho,t) \,d\rho \quad \text{for } s \in [0,R^n] \text{ and } t \in [0,\TM),
\end{equation}
which implies that 
\begin{equation}\label{FormuleW}
U_s(s,t) = \frac{1}{n} u(s^{\frac{1}{n}},t) \quad \text{and} \quad U_{ss}(s,t)=\frac{1}{n^2} s^{\frac{1}{n}-1} u_r(s^{\frac{1}{n}},t) \quad \text{for } s \in (0,R^n) \text{ and } t \in (0,\TM).
\end{equation}
By the definition of $U$ and by exploiting \eqref{FormuleW}, we obtain 
\begin{equation}\label{Stima_U_t}
\begin{split}
U_t(s,t) &= s^{1-\frac{1}{n}}(u(s^{\frac{1}{n}},t)+1)^{m_1-1} u_r(s^{\frac{1}{n}},t)
- s^{1-\frac{1}{n}} u(s^{\frac{1}{n}},t) (u(s^{\frac{1}{n}},t)+1)^{m_2-1} z_r(s^{\frac{1}{n}},t)\\
&+\lambda \int_0^{s^{\frac{1}{n}}} \rho^{n-1} u(\rho,t) \,d\rho
- \mu \int_0^{s^{\frac{1}{n}}} \rho^{n-1} u^k(\rho,t) \,d\rho\\
&= n^2 s^{2-\frac{2}{n}}(n U_s(s, t)+1)^{m_1-1} U_{ss}(s,t)
- n s^{1-\frac{1}{n}} U_s(s,t) (n U_s(s, t)+1)^{m_2-1} z_r(s^{\frac{1}{n}},t)\\ 
&+ \lambda U(s,t) - \mu n^{k-1} \int_0^s U_s^k(\sigma,t)\,d\sigma \quad \text{for } s \in (0,R^n) \text{ and } t \in (0,\TM).
\end{split}
\end{equation}
Besides, with an integration over $(0,r)$ of the second equation in \eqref{RadialProblem}  and the substitution $r=s^{\frac{1}{n}}$, we arrive at  
\[
z_r= \frac{m(t)}{n} s^{\frac{1}{n}} - \frac{1}{n} s^{\frac{1}{n}-1} \int_0^s f(n U_s(\sigma,t))\,d\sigma \quad \text{for } s \in (0,R^n) \text{ and } t \in (0,\TM),
\]
which inserted into relation \eqref{Stima_U_t} gives (observe $U\geq 0$)
\begin{equation}\label{StimaDisUt}
\begin{split}
U_t(s,t) &\geq n^2 s^{2-\frac{2}{n}}(n U_s+1)^{m_1-1} U_{ss} 
- s U_s (n U_s+1)^{m_2-1} m(t) + U_s (n U_s+1)^{m_2-1} \int_0^s f(n U_s(\sigma,t))\,d\sigma\\
&\quad - \mu n^{k-1} \int_0^s U_s^k(\sigma,t)\,d\sigma \quad \text{for all } s \in (0,R^n) \text{ and } t \in (0,\TM).
\end{split}
\end{equation}
In addition, given $s_0 \in (0, R^n)$, $\gamma \in (-\infty, 1)$ and $U$ as in \eqref{W}, we introduce the moment-type functional
\begin{equation}\label{Phi}
\phi(t):=\int_0^{s_0} s^{-\gamma} (s_0-s) U(s,t)\,ds \quad \text{for } t \in [0,\TM),
\end{equation}
which is well defined and belongs to $C^0([0,\TM)) \cap C^1((0,\TM))$.
Moreover, we define 
\begin{equation}\label{psi}
\psi(t):= \int_0^{s_0} s^{1-\gamma} (s_0-s) U_s^{m_2+\alpha}(s,t) \,ds  \quad \text{for } t \in (0,\TM),
\end{equation}
and the set
\begin{equation}\label{S_phi}
S_{\phi}=\left\{t \in (0,\TM): \phi(t) \geq \frac{M-s_0}{(1-\gamma)(2-\gamma) \omega_n}s_0^{2-\gamma}\right\},
\end{equation}
where $M$ is the bound of the $L^1(\Omega)$-norm of $u$ established in \eqref{boundednessMass} and $
\omega_n=n |B_1(0)|$.
With these preparations in our hands, let us give a series of necessary lemmas, some of which are not new.

We start with a result dealing with the concavity of $U$ and some estimate for $m(t)$.
\resetconstants{c}
\begin{lemma}\label{Concav}
Let $f_1, f_2$ and $u_0$ satisfy \eqref{f_i_u0_GENERALE} and  \eqref{f_iDef}, $\alpha>\beta$ and $\gamma\in \left(-\infty,1\right).$ Then the following relations hold:
\[
U_{ss}(s,t) \leq 0 \quad \text{for all } s \in (0,R^n) \text{ and } t \in (0, \TM),
\]
and
\begin{equation}\label{m(t)}
m(t) \leq c_0 + \frac{1}{2s}  \int_0^s f(n U_s(\sigma,t))\,d\sigma \quad \text{for }
s_0 \in \left(0,\frac{R^n}{6}\right] \text{and for all } s \in (0,s_0),
\end{equation}
for all $t \in S_{\phi}$, and with some $c_0$.
\begin{proof}
As to the concavity property, the proof is based on minor adjustments of \cite[Lemma 3.2]{wang2023blow} or \cite[Lemma 2.2]{WinklerNoNLinearanalysisSublinearProduction}. The remaining conclusion follows by closely reasoning as in  \cite[Lemma 2.5]{LiuLiBlowUpAttr-Rep};  we herein only add that $\alpha>\beta$ is required just for the construction of 
 $$c_0= \chi f_1\left(\frac{8n}{2^{\gamma} (3-\gamma)\omega_n}\right)+\frac{1}{6}
\left(\frac{\chi k_3 (\alpha-\beta)}{\beta} \left(\frac{2\xi k_2 \beta}{\chi k_3 \alpha}
\right)^{\frac{\alpha}{\alpha-\beta}} + L (\chi + 2)\right)>0,$$ being $L=L(K)$ estimated by \begin{equation*}
|f(s)|, f_1(s), f_2(s) \leq L \quad \text{for } s \in [0,K], K>0.
\end{equation*}
(We point out that $c_0$ will be used in some other places below.)
\end{proof}
\end{lemma}
\resetconstants{c}
Let us now start with the analysis of the temporal evolution of $\phi$.
\resetconstants{c}
\begin{lemma}\label{Phi_T}
Under the same assumptions of Lemma \ref{Concav}, let $s_0 \in \left(0,\frac{R^n}{6}\right]$. Then 
\begin{equation}\label{Phi1_T}
\begin{split}
\phi'(t) &\geq n^2 \int_0^{s_0} s^{2-\frac{2}{n}-\gamma} (s_0-s) (n U_s+1)^{m_1-1} U_{ss} \,ds
-c_0 \int_0^{s_0} s^{1-\gamma} (s_0-s) U_s (n U_s+1)^{m_2-1}\,ds\\
&\quad +\frac{1}{2} \int_0^{s_0} s^{-\gamma} (s_0-s) U_s (n U_s+1)^{m_2-1} 
\left[\int_0^s f(n U_s(\sigma,t))\,d\sigma\right]\,ds\\
&\quad - \mu n^{k-1} \int_0^{s_0} s^{-\gamma} (s_0-s) \left[ \int_0^s U_s^k(\sigma,t)\,d\sigma \right]\, ds \\
&= I_1 + I_2 + I_3 + I_4 \quad \text{for all }  t \in S_{\phi}.
\end{split}
\end{equation}
\begin{proof}
By the definition of $\phi$ (recall \eqref{Phi}) and exploiting \eqref{StimaDisUt}, we get this estimate 
\begin{equation*}
\begin{split}
\phi'(t)&= \int_0^{s_0} s^{-\gamma} (s_0-s) U_t(s,t)\,ds\\ 
&\geq n^2 \int_0^{s_0} s^{2-\frac{2}{n}-\gamma} (s_0-s) (n U_s+1)^{m_1-1} U_{ss} \,ds
- \int_0^{s_0} s^{1-\gamma} (s_0-s) U_s (n U_s+1)^{m_2-1}m(t)\,ds\\
&\quad +\int_0^{s_0} s^{-\gamma} (s_0-s) U_s (n U_s+1)^{m_2-1} 
\left[\int_0^s f(n U_s(\sigma,t))\,d\sigma\right]\,ds\\
&\quad - \mu n^{k-1} \int_0^{s_0} s^{-\gamma} (s_0-s) \left[ \int_0^s U_s^k(\sigma,t)\,d\sigma \right]\, ds \quad \text{for all }  t \in S_{\phi}.
\end{split}
\end{equation*}
Finally, by applying \eqref{m(t)}, we obtain the thesis.
\end{proof}
\end{lemma}
The next results provide some lower bounds of $I_1,I_2,I_3$ and $I_4$ in respect of $\psi(t)$ defined in \eqref{psi}.
\resetconstants{c}
\begin{lemma}\label{I1-I4}
Let $u_0$ satisfy the related assumptions in \eqref{f_i_u0_GENERALE} and \eqref{f_iDef},  $m_1 \in \R$, $m_2, \alpha>0$, $k>1$ and $\gamma \in (-\infty,1)$.
\begin{itemize}[label={\ding{226}}]
\item If $m_1, m_2,  \alpha$ and $\gamma$ comply with  
\[
m_2+\alpha >m_1 \quad \text{and} \quad \gamma > 2- \frac{2}{n} \frac{(m_2+\alpha)}{(m_2+\alpha-m_1)} \quad \text{if} \quad m_1 \geq 0,
\]
\[
\qquad \qquad \qquad \qquad \qquad \qquad \qquad \qquad \gamma < 2- \frac{2}{n} \quad \text{if} \quad m_1 < 0,
\]
then there exist $\varepsilon>0$ (sufficiently small), and $\const{Sil3}$, $\const{Sil4}$, $\const{Sil5}$, $\const{Sil6}$, $\const{Sil7}$ such that for any $s_0 \in \left(0,\frac{R^n}{6}\right]$ 
\[
I_1 \geq
\begin{cases}
-\const{Sil3} s_0^{(3-\gamma)\frac{(m_2+\alpha-m_1)}{m_2+\alpha}-\frac{2}{n}}  
\psi^{\frac{m_1}{m_2+\alpha}}(t) - \const{Sil4} s_0^{3-\gamma-\frac{2}{n}}  & \text{if } m_1 >0,\\
-\const{Sil5} s_0^{(3-\gamma)\frac{(m_2+\alpha-\varepsilon)}{m_2+\alpha}-\frac{2}{n}}  
\psi^{\frac{\varepsilon}{m_2+\alpha}}(t) - \const{Sil6} s_0^{3-\gamma-\frac{2}{n}}  & \text{if } m_1 =0,\\
- \const{Sil7} s_0^{3-\gamma-\frac{2}{n}}  & \text{if } m_1 <0
\end{cases}
\]
for all $t \in S_{\phi}$.
\item If $m_2, \alpha, k$ and $\gamma$ are such that   
 \[
m_2+\alpha >k \quad \text{and} \quad 2-\frac{(m_2+\alpha)}{k}<\gamma<1,
\]
then there exists $\const{Sil8}$ such that for any $s_0 \in \left(0,\frac{R^n}{6}\right]$ 
\[
I_4 \geq -\const{Sil8} \psi^{\frac{k}{m_2+\alpha}}(t) s_0^{(3-\gamma)\frac{(m_2+\alpha-k)}{m_2+\alpha}} \quad \text{for all } t \in S_{\phi}.
\]
\end{itemize}
\begin{proof}
The proof can be found in \cite[Lemmas 3.6 and 3.7]{TANAKA2021}.
\end{proof}
\end{lemma}
As to the estimate of $I_2+I_3$ we need to rearrange some computations; this is exactly where we have to go beyond the analysis in \cite{TANAKA2021} and \cite{wang2023blow}.
\resetconstants{c}
\begin{lemma}\label{I2-I3}
Under the same assumptions of Lemma \ref{Concav}, let moreover $m_2, \alpha>0$ comply with $m_2+\alpha>1.$
Then for some $\const{Sil11}$, $\const{Sil12}$ and $\const{Sil13}$  we have 
\begin{equation}\label{I2I3}
I_2+I_3 \geq \const{Sil11} \psi(t) -\const{Sil12} \psi^{\frac{m_2}{m_2+\alpha}}(t) 
s_0^{(3-\gamma)\frac{\alpha}{m_2+\alpha}} - \const{Sil13} s_0^{3-\gamma} \; \textrm{ for any } s_0 \in \left(0,\frac{R^n}{6}\right] \textrm{ and for all } t \in S_{\phi}.
\end{equation}

\begin{proof}
Since $\alpha>\beta$, by applying the Young inequality and \eqref{f_iDef}, we get 
\[
\xi f_2(u) \leq \xi k_2 (u+1)^{\beta} \leq \frac{\chi k_3}{2} (u+1)^{\alpha} + \const{Sil9}
\leq \frac{\chi}{2} f_1(u) + \const{Sil9},
\]
with $\const{Sil9}=\left(\frac{2\beta \xi k_2}{\chi k_3 \alpha}\right)^{\frac{\alpha}{\alpha-\beta}}\frac{\chi k_3 (\alpha-\beta)}{2 \beta}>0$, which implies
\begin{equation}\label{F}
\frac{\chi}{2} f_1(u) - \const{Sil9} \leq f(u) \leq \chi f_1(u).
\end{equation}
From the concavity of $U$ in Lemma \ref{Concav}, it is seen that $U_s$ is nonincreasing, namely $U_s(\sigma, t) \geq U_s(s,t)$ for any $\sigma \in (0, s)$.
Henceforth, since  $f_1$ is nondecreasing  (recall \eqref{f_iDef}),  we have
\begin{equation}\label{Int_f}
\int_0^s f_1(n U_s(\sigma,t))\,d\sigma \geq \int_0^s f_1(n U_s(s,t))\,d\sigma
= s f_1(n U_s(s,t)) \quad \text{for all } s \in (0,s_0)  \text{ and } t \in (0,\TM).
\end{equation}
Therefore, from \eqref{F} and \eqref{Int_f} we derive the estimate 
\[
\begin{split}
I_3 &= \frac{1}{2} \int_0^{s_0} s^{-\gamma} (s_0-s) U_s (n U_s+1)^{m_2-1} 
\left[\int_0^s f(n U_s(\sigma,t))\,d\sigma\right]\,ds\\
&\geq  \frac{1}{2} \int_0^{s_0} s^{-\gamma} (s_0-s) U_s (n U_s+1)^{m_2-1} 
\left(\frac{\chi}{2} s f_1(nU_s)-\const{Sil9} s \right)\,ds\\
&\geq \frac{\chi k_3}{4} \int_0^{s_0} s^{1-\gamma} (s_0-s) U_s (n U_s+1)^{m_2-1+\alpha}\,ds - \frac{\const{Sil9}}{2} \int_0^{s_0} s^{1-\gamma} (s_0-s) U_s (n U_s+1)^{m_2-1}\,ds \quad \text{for all }  t \in S_{\phi}.
\end{split}
\]
Then, for $c_0$ as in Lemma \ref{Concav},
\begin{equation}\label{Sum2_3}
I_2 + I_3 \geq \frac{\chi k_3}{4} \int_0^{s_0} s^{1-\gamma} (s_0-s) U_s (n U_s+1)^{m_2-1+\alpha}\,ds - \left(c_0 + \frac{\const{Sil9}}{2}\right) \int_0^{s_0} s^{1-\gamma} (s_0-s) U_s (n U_s+1)^{m_2-1}\,ds \quad \text{for all }  t \in S_{\phi}.
\end{equation}
Now, we focus on the first integral in \eqref{Sum2_3}. Since $m_2+\alpha-1>0$, clearly we get that $(n U_s + 1)^{m_2 + \alpha -1} > (n U_s)^{m_2+\alpha-1}$, and by exploiting this latter we have 
\begin{equation}\label{1T}
\frac{\chi k_3}{4} \int_0^{s_0} s^{1-\gamma} (s_0-s) U_s (n U_s+1)^{m_2-1+\alpha}\,ds \geq \frac{\chi k_3 n^{m_2+\alpha-1}}{4}  \int_0^{s_0} s^{1-\gamma} (s_0-s)  U_s^{m_2+\alpha}\,ds = \frac{\chi k_3 n^{m_2+\alpha-1}}{4} \psi(t) \quad \text{on }   S_{\phi}.
\end{equation}
In order to treat the second integral in \eqref{Sum2_3},  we use inequality \eqref{algebricIne}: since $m_2>0$, we obtain for $\const{Sil103}=\left(c_0 + \frac{\const{Sil9}}{2}\right)$ the following estimate
\begin{equation}\label{2T}
\begin{split}
&- \const{Sil103} \int_0^{s_0} s^{1-\gamma} (s_0-s) U_s (n U_s+1)^{m_2-1}\,ds
= - \frac{1}{n} \const{Sil103}\int_0^{s_0} s^{1-\gamma} (s_0-s) \frac{n U_s}{n U_s+1} (n U_s+1)^{m_2}\,ds\\
& \geq - \frac{1}{n} \const{Sil103} \int_0^{s_0} s^{1-\gamma} (s_0-s)(n 
U_s+1)^{m_2}\,ds\\
& \geq - \const{Sil14} \int_0^{s_0} s^{1-\gamma} (s_0-s) U_s^{m_2}\,ds
-  \const{Sil15} \int_0^{s_0} s^{1-\gamma} (s_0-s)\,ds \quad \text{for all }  t \in S_{\phi}.
\end{split}
\end{equation}
Additionally, from $\frac{m_2}{m_2+\alpha}<1$ we can apply the H\"{o}lder inequality to the first integral in \eqref{2T}, which leads to
\begin{equation}\label{3T}
\begin{split}
\int_0^{s_0} s^{1-\gamma} (s_0-s) U_s^{m_2}\,ds
&= \int_0^{s_0} \left(s^{1-\gamma} (s_0-s) U_s^{m_2+\alpha}\right)^{\frac{m_2}{m_2+\alpha}} s^{\frac{(1-\gamma)\alpha}{m_2+\alpha}} (s_0-s)^{\frac{\alpha}{m_2+\alpha}}\,ds\\
&\leq \left(\int_0^{s_0} s^{1-\gamma} (s_0-s) U_s^{m_2+\alpha}\,ds\right)^{\frac{m_2}{m_2+\alpha}} \left(\int_0^{s_0} s^{1-\gamma} (s_0-s)\,ds\right)^{\frac{\alpha}{m_2+\alpha}}\\
&=\psi^{\frac{m_2}{m_2+\alpha}}(t) \, ((2-\gamma)(3-\gamma))^{-\frac{\alpha}{m_2+\alpha}} s_0^{\frac{(3-\gamma)\alpha}{m_2+\alpha}} \quad \text{for all }  t \in S_{\phi}.
\end{split}
\end{equation}
By putting \eqref{3T} into \eqref{2T}, we obtain for every $t\in S_{\phi}$
\begin{equation}\label{4T}
\begin{split}
- \const{Sil103} \int_0^{s_0} s^{1-\gamma} (s_0-s) U_s (n U_s+1)^{m_2-1}\,ds
&\geq - \const{Sil14} \psi^{\frac{m_2}{m_2+\alpha}}(t) \, ((2-\gamma)(3-\gamma))^{-\frac{\alpha}{m_2+\alpha}} s_0^{\frac{(3-\gamma)\alpha}{m_2+\alpha}} - \frac{\const{Sil15}}{(2-\gamma)(3-\gamma)}  s_0^{3-\gamma},
\end{split}
\end{equation}
so by invoking \eqref{Sum2_3} and taking into account \eqref{1T} and \eqref{4T}, we can conclude. 
\end{proof}
\end{lemma} 
In order to obtain the desired superlinear ODI for $\phi$, we have to rely on some relations involving  $U$ and $\psi$ and $\phi.$
\resetconstants{c}
\begin{lemma}\label{Phi-Psi}
Let $u_0$ satisfy its related assumptions in  \eqref{f_i_u0_GENERALE} and \eqref{f_iDef}, $m_2, \alpha>0$ and $\gamma \in (-\infty,1)$ be such that 
\[
m_2+\alpha >1 \quad \text{and} \quad 2-(m_2+\alpha)<\gamma<1.
\]
Then there exist $\const{Sil16}$, $\const{Sil17}$ such that for any $s_0 \in \left(0,\frac{R^n}{6}\right]$ 
\begin{equation}\label{Phi-Psi1}
U(s,t) \leq \const{Sil16} s^{\frac{m_2+\alpha+\gamma-2}{m_2+\alpha}} (s_0-s)^{-\frac{1}{m_2+\alpha}} \psi^{\frac{1}{m_2+\alpha}}(t) \quad \text{for all } s \in (0,s_0) \text{ and }  
t \in S_{\phi},
\end{equation}
and
\begin{equation}\label{Phi-Psi2}
\psi(t) \geq \const{Sil17} s_0^{-(3-\gamma)(m_2+\alpha-1)} \phi^{m_2+\alpha}(t) \quad 
\text{for all } t \in S_{\phi}.
\end{equation}
\begin{proof}
The proof of \eqref{Phi-Psi1} and \eqref{Phi-Psi2} follows by 
\cite[Lemma 3.8]{TANAKA2021} and \cite[Lemma 3.7]{wang2021boundedness}, respectively.
\end{proof}
\end{lemma}
\resetconstants{c}
The following is precisely the lemma relying on assumption \ref{BU}.
\begin{lemma}\label{ODI}
Under the same assumptions of Lemma \ref{Concav}, let $m_1 \in \R$, $m_2, \alpha, \beta>0$ and $k>1$ be such that constrains \ref{BU} in Assumptions \ref{generalassumptions} are satisfied. 
Then there exist $\varepsilon>0$ small enough,
$\gamma=\gamma(m_1, m_2, \alpha, k) \in (-\infty,1)$ and $\const{Sil20}$, $\const{Sil21}$, $\const{Sil22}$, $\const{Sil23}$,  $\const{Sil24}$ and $\const{Sil25}$
 such that for $s_0 \in \left(0,\frac{R^n}{6}\right]$ one has  
\[
\phi'(t) \geq
\begin{cases}
\const{Sil20} s_0^{-(3-\gamma)(m_2+\alpha-1)} \phi^{m_2+\alpha}(t) - \const{Sil21} 
s_0^{3-\gamma-\frac{2}{n}\frac{(m_2+\alpha)}{m_2+\alpha-m_1}}  & \text{if } m_1 >0,\\
\const{Sil22} s_0^{-(3-\gamma)(m_2+\alpha-1)} \phi^{m_2+\alpha}(t) - \const{Sil23}
s_0^{3-\gamma-\frac{2}{n}\frac{(m_2+\alpha)}{m_2+\alpha-\epsilon}}  & \text{if } m_1 =0,\\
\const{Sil24} s_0^{-(3-\gamma)(m_2+\alpha-1)} \phi^{m_2+\alpha}(t) - \const{Sil25}
s_0^{3-\gamma-\frac{2}{n}}  & \text{if } m_1 < 0,
\end{cases}
\quad \quad \textrm{for all } t \in S_{\phi}.
\]

\begin{proof}
By substituting the estimates of $I_1, I_2, I_3, I_4$ given in Lemmas \ref{I1-I4}, \ref{I2-I3} and \ref{Phi-Psi} into Lemma \ref{Phi_T},  we can adapt \cite[Lemma 3.10]{TANAKA2021} taking into account that \cite[(3.12), Lemma 3.5]{TANAKA2021} is replaced by \eqref{I2I3}, so being necessary manipulating the term involving $\psi^\frac{m_2}{m_2+\alpha}(t)s_0^{(3-\gamma)\frac{\alpha}{m_2+\alpha}}$ by the Young inequality and relation \eqref{Phi-Psi2}.
\end{proof}
\end{lemma}
The previous lemmata allows us to conclude. 
\subsubsection*{Proof of Theorem \ref{maintheorem2}:}
We focus only on the situation where $m_1>0$, the cases $m_1=0$ and $m_1<0$ being similar. Since \ref{BU} holds, we can apply Lemma \ref{ODI} and find $\gamma \in (-\infty,1)$, $\const{Sil30}, \const{Sil31}$ such that for each $u_0$ fulfilling the related restrictions in \eqref{f_i_u0_GENERALE} and \eqref{f_iDef}, and $s_0 \leq \frac{R^n}{6}$, one deduces
\begin{equation}\label{ODI1}
\phi'(t) \geq \const{Sil30} s_0^{-(3-\gamma)(m_2+\alpha-1)} \phi^{m_2+\alpha}(t) - \const{Sil31} 
s_0^{3-\gamma-\frac{2}{n}\frac{(m_2+\alpha)}{m_2+\alpha-m_1}}  \quad \text{for all } 
t \in S_{\phi}.
\end{equation}
Next we pick $s_0 \leq \frac{R^n}{6}$ small enough such that
\begin{equation}\label{Cond1_s0}
s_0\leq \frac{M}{2} 
\end{equation}
and
\begin{equation}\label{Cond2_s0}
s_0^{m_2+\alpha-\frac{2}{n}\frac{(m_2+\alpha)}{m_2+\alpha-m_1}} \leq 
\frac{\const{Sil30}}{2 \const{Sil31}} \left(\frac{M}{2(1-\gamma)(2-\gamma)\omega_n}\right)^{m_2+\alpha}.
\end{equation}
Moreover, from
\[
\frac{M_0-\epsilon_0}{\omega_n} \int_{s_*}^{s_0} s^{-\gamma} (s_0-s)\,ds=\frac{M_0-\epsilon_0}{\omega_n}\left[\frac{s_0^{2-\gamma}}{(1-\gamma)(2-\gamma)}-s_*^{1-\gamma}\frac{s_0}{1-\gamma}+\frac{s_*^{2-\gamma}}{2-\gamma}\right],
\]
we can take $\epsilon_0\in \left(0,\frac{s_0}{2}\right)$  and $s_*\in \left(0,s_0\right)$ so small and  satisfying  
\begin{equation}\label{CondInt}
\frac{M_0-\epsilon_0}{\omega_n} \int_{s_*}^{s_0} s^{-\gamma} (s_0-s)\,ds > 
\frac{M_0-s_0}{(1-\gamma)(2-\gamma)\omega_n} s_0^{2-\gamma}.
\end{equation}
We set $r_*:= s_*^{\frac{1}{n}} \in (0,R)$ and require to $u_0$ to satisfy \eqref{M0}. In order to show that $\TM <\infty$, we argue by establishing that if by contradiction  $\TM=\infty$, then $\tilde{T}=\sup \tilde{S}\in (0,\infty]$, with  
\begin{equation}\label{SetS}
\tilde{S}:= \left\{T>0 \textrm{ such that } \phi(t) > \frac{M-s_0}{(1-\gamma)(2-\gamma)\omega_n} s_0^{2-\gamma} \quad \text{for all } t \in [0,T] \right\},
\end{equation}
would be at the same time finite and infinite. 

First, by observing that $\phi(0) > \frac{M-s_0}{(1-\gamma)(2-\gamma)\omega_n} s_0^{2-\gamma}$, the continuity of $\phi$ ensures that $\tilde{S}$ is not empty. Indeed, we have that  for any $s \in (s_*, R^n)$
\[
U(s,0) \geq U(s_*,0)= \frac{1}{\omega_n} \int_{B_{r_*}(0)} u_0\,dx \geq  \frac{M_0-\epsilon_0}{\omega_n}.
\]
In turn, due to $M=\max\{M_0, C\}=M_0$ thanks to the assumption $M_0>C$, we deduce from \eqref{CondInt} that 
\[
\begin{split}
\phi(0)&\geq \int_{s_*}^{s_0} s^{-\gamma} (s_0-s) U(s,0)\,ds
\geq \frac{M_0-\epsilon_0}{\omega_n} \int_{s_*}^{s_0} s^{-\gamma} (s_0-s)\,ds\\
& > \frac{M_0-s_0}{(1-\gamma)(2-\gamma)\omega_n} s_0^{2-\gamma}
= \frac{M-s_0}{(1-\gamma)(2-\gamma)\omega_n} s_0^{2-\gamma}.
\end{split}
\]
Now, by exploiting \eqref{SetS} and \eqref{Cond1_s0} leads to
\[
\phi(t) \geq \frac{M}{2(1-\gamma)(2-\gamma)\omega_n} s_0^{2-\gamma} \quad \text{for all } t \in (0,\tilde{T}),
\]
so that $(0,\tilde{T}) \subset S_{\phi}$ (recall \eqref{S_phi}). 

Secondly, condition \eqref{Cond2_s0} implies 
\[
\frac{\frac{\const{Sil30}}{2} s_0^{-(3-\gamma)(m_2+\alpha-1)} \phi^{m_2+\alpha}(t)}
{\const{Sil31} s_0^{3-\gamma-\frac{2}{n}\frac{(m_2+\alpha)}{m_2+\alpha-m_1}}}
\geq \frac{\const{Sil30}}{2 \const{Sil31}} \left(\frac{M}{2(1-\gamma)(2-\gamma)\omega_n}\right)^{m_2+\alpha} s_0^{-(m_2+\alpha)+\frac{2}{n}\frac{(m_2+\alpha)}{m_2+\alpha-m_1}} \geq 1 \quad \textrm{on} \; (0,\tilde{T}),
\]
which subsequently gives from \eqref{ODI1} 
\begin{equation}\label{ODI2}
\phi'(t) \geq \frac{\const{Sil30}}{2} s_0^{-(3-\gamma)(m_2+\alpha-1)} \phi^{m_2+\alpha}(t) \geq 0 \quad \text{for all } t \in (0, \tilde{T}).
\end{equation}
With these facts, let us now establish the inconsistency $\tilde{T}<\infty$ and $\tilde{T}=\infty$. By integrating inequality \eqref{ODI2} on $(0,\tilde{T})$, we have that 
\[
\int_0^{\tilde{T}} \left(\frac{1}{1-(m_2+\alpha)} \phi^{1-(m_2+\alpha)}(t)\right)'\,dt \geq 
\int_0^{\tilde{T}} \frac{\const{Sil30}}{2} s_0^{-(3-\gamma)(m_2+\alpha-1)} \, dt,
\]
so that due to $m_2+\alpha-1>0$ and the nonnegative property of $\phi$
\[
\frac{\const{Sil30}}{2} s_0^{-(3-\gamma)(m_2+\alpha-1)} \tilde{T}
\leq \frac{\phi^{1-(m_2+\alpha)}(\tilde{T})}{1-(m_2+\alpha)} - 
\frac{\phi^{1-(m_2+\alpha)}(0)}{1-(m_2+\alpha)} 
\leq \frac{\phi^{1-(m_2+\alpha)}(0)}{m_2+\alpha-1}, 
\]
or explicitly 
\[
\tilde{T} \leq \frac{2}{\const{Sil30} (m_2+\alpha-1) \phi^{m_2+\alpha-1}(0)} 
s_0^{(3-\gamma)(m_2+\alpha-1)} < \infty.
\]
Nevertheless, we have to exclude the finiteness of $\tilde{T}=\sup \tilde{S}$; in fact, by the definition of $\tilde{S}$ in \eqref{SetS}, we should have 
\begin{equation}\label{SupSTilde}
\phi(\tilde{T})= \frac{M-s_0}{(1-\gamma)(2-\gamma)\omega_n} s_0^{2-\gamma},
\end{equation}
because if 
\[
\phi(\tilde{T})> \frac{M-s_0}{(1-\gamma)(2-\gamma)\omega_n} s_0^{2-\gamma},
\]
by continuity of $\phi$  we would have that $\tilde{T}$ cannot be the supremum of $\tilde{S}$. But \eqref{SupSTilde} cannot be true since from the nondecreasing of $\phi$ in view of inequality \eqref{ODI2}, we would arrive at
\[
\frac{M-s_0}{(1-\gamma)(2-\gamma)\omega_n} s_0^{2-\gamma}=\phi(\tilde{T})\geq \phi(0)>\frac{M-s_0}{(1-\gamma)(2-\gamma)\omega_n} s_0^{2-\gamma}. 
\]
As a conclusion, the constructed $u_0$ implies that $\TM$ has to be finite.
\qed
\begin{remark}[Finite-time blow-up with $m_2 \neq m_3$]\label{Casem2neqm3}
As already said in $\S$\ref{section3}, the proof of the blow-up for problem \eqref{problem3x3} with $m_2 \neq m_3$ is still an open problem. Since the transformation $z=\chi v -\xi w$ used above to reorganize problem \eqref{problem3x3} with $m_2=m_3$ into the  simplified version \eqref{ReducedProblem} is not longer employable,  by reasoning as in Lemma \ref{Phi_T}, the corresponding inequality \eqref{Phi1_T} would read 
\[
\begin{split}
\phi'(t) &\geq n^2 \int_0^{s_0} s^{2-\frac{2}{n}-\gamma} (s_0-s) (n U_s+1)^{m_1-1} U_{ss} \,ds
-\chi f_{1_{\gamma}} \int_0^{s_0} s^{1-\gamma} (s_0-s) U_s (n U_s+1)^{m_2-1}\,ds\\
&\quad +\frac{\chi}{2} \int_0^{s_0} s^{-\gamma} (s_0-s) U_s (n U_s+1)^{m_2-1} 
\left[\int_0^s f_1(n U_s(\sigma,t))\,d\sigma\right]\,ds
- \mu n^{k-1} \int_0^{s_0} s^{-\gamma} (s_0-s) \left[ \int_0^s U_s^k(\sigma,t)\,d\sigma \right]\, ds \\
& \quad +\xi \int_0^{s_0} s^{1-\gamma} (s_0-s) U_s (n U_s+1)^{m_3-1} m_2(t)\,ds
-\xi \int_0^{s_0} s^{-\gamma} (s_0-s) U_s (n U_s+1)^{m_3-1} 
\left[\int_0^s f_2(n U_s(\sigma,t))\,d\sigma\right]\,ds,
 \end{split}
\]
 valid for  all 
$t \in S_{\phi}$ and with  $f_{1_{\gamma}}= f_1\left(\frac{8n}{2^{\gamma} (3-\gamma)\omega_n}\right)$. In particular, the extra terms involving the repulsion coefficient $\xi$ make the analysis more complex. This is also connected to the signs of such terms (exactly opposite than those of the integrals associated to the attraction coefficient $\chi$), which do not allow to use the right inequalities tied to the hypotheses on $f_i$, and in turn on $m_i$. 
\end{remark}
\subsubsection*{Acknowledgements}
The authors are members of the Gruppo Nazionale per l'Analisi Matematica, la Probabilit\`a e le loro Applicazioni (GNAMPA) of the Istituto Nazionale di Alta Matematica (INdAM), and are partially supported by the research projects \emph{Evolutive and Stationary Partial Differential Equations with a Focus on Biomathematics} (2019, Grant Number: F72F20000200007), {\em  Analysis of PDEs in connection with real phenomena} (2021, Grant Number: F73C22001130007), funded by  \href{https://www.fondazionedisardegna.it/}{Fondazione di Sardegna}. GV is also supported by MIUR (Italian Ministry of Education, University and Research) Prin 2017 \emph{Nonlinear Differential Problems via Variational, Topological and Set-valued Methods} (Grant Number: 2017AYM8XW). 

\end{document}